\documentclass[10pt]{amsart}
\usepackage{latexsym}
\usepackage{amsmath}
\usepackage{amssymb,amsthm}
\usepackage[french]{babel}
\usepackage{comment}
\usepackage{graphicx,epic,eepic}
\usepackage{mudessin}

\newtheorem{theo}{Th\'eor\`eme}
\newtheorem{lem}{Lemme}[section]

\newtheorem{corol}[lem]{Corollaire}
\newtheorem{prop}[lem]{Proposition}
\theoremstyle{definition}
\newtheorem{remarque}[lem]{Remarque}

\def\RR{\mathbb{R}}
\def\NN{\mathbb{N}}

\def\eps{\varepsilon}
\def\DDD{\mathcal{D}}
\def\CCC{\mathcal{C}}
\def\III{\mathcal{I}}
\def\PPP{\mathcal{P}}

\def\SSS{\mathcal{S}}
\def\pp{\mathfrak{p}}
\def\gg{\mathfrak{g}}

\def\uk{\mathsf{u}_k}
\def\fk{\mathsf{f}_k}
\def\gk{\mathsf{g}_k}
\def\vk{\mathsf{v}_k}
\def\wk{\mathsf{w}_k}

\def\fkbar{\overline{\mathsf{f}}_k}

\def\vkbar{\overline{\mathsf{v}}_k}
\def\wkbar{\overline{\mathsf{w}}_k}
\def\ukn{\mathsf{u}_{kn}}
\def\fkn{\mathsf{f}_{kn}}
\def\gkn{\mathsf{g}_{kn}}
\def\vkn{\mathsf{v}_{kn}}
\def\wkn{\mathsf{w}_{kn}}

\def\uknu{\mathsf{u}_{kn}^{\star}}
\def\fknu{\mathsf{f}_{kn}^{\star}}

\def\nuu{\nu^{\star}}
\def\estar{e^{\star}}

\def\indic{1\!\!1}
\def\d{\,{ \rm d}}
\def\eps{\varepsilon}

\def\surpid{\frac{1}{(2\pi)^d}}

\def\egaldef{:=}
\def\defegal{=:}

\def\ubar{\overline{u}}

\def\Re{\text{Re}\,}
\def\Im{\text{Im}\,}

\def\Ld{L^2(\RR^d)}
\def\Cio{C^\infty_{0}}
\def\loc{\rm loc}

\newcommand{\abs}[1]{|#1|}
\newcommand{\carre}[1]{\abs{#1}^2}
\newcommand{\norme}[1]{\|#1\|}

\def\supp{\text{supp}\,}

\title{Op\'erateur de Schr\"odinger avec potentiel singulier multipolaire}
\author{Thomas Duyckaerts}
\date{\today}

\begin{document}
\maketitle
\begin{abstract}
On \'etudie un op\'erateur de la forme: $ -\Delta+V$ sur $\RR^d$, o\`u $V$ est un potentiel admettant plusieurs p\^oles en $a/r^2$. Plus pr\'ecis\'ement, on d\'emontre l'estimation de r\'esolvante tronqu\'ee \`a hautes fr\'equences, classique dans les cas non-captifs, et qui implique des propri\'et\'es de r\'egularisation sur l'\'equation de Schr\"odinger correspondante. La preuve est bas\'ee sur l'introduction d'une mesure de d\'efaut micro-locale semi-classique. 
\end{abstract}

\section{Introduction}
Nous consid\'erons dans ce travail un op\'erateur de la forme:
\begin{equation}
P=-\Delta+V,
\end{equation}
sur l'espace $\RR^d$ ($d\geq2$), o\`u:  
$$ \Delta=\sum_j\frac{\partial^2}{\partial_{x_j}^2}$$
est le Laplacien standard, et le potentiel $V$ est une fonction r\'eelle d\'efinie sur $\RR^d$. L'\'etude de la r\'esolvante:
$$ R_z=(P-z)^{-1},\; z\notin \RR,$$
pr\`es de l'axe r\'eel, int\'eressante en elle m\^eme, permet aussi de pr\'eciser le comportement des solutions des \'equations d'onde et de Schr\"odinger associ\'ees \`a $P$. Lorsque $V$ est r\'egulier, de nombreuses in\'egalit\'es ont \'et\'e d\'emontr\'ees sur des normes de $R_z$ dans des espaces \`a poids, notamment l'in\'egalit\'e standard \`a haute fr\'equence, sur la r\'esolvante tronqu\'ee:
\begin{equation}
\label{ineg HF standard}
|| \chi R_z \chi||_{L^2\rightarrow L^2} \leq \frac{C}{1+|z|}
\end{equation}
o\`u $z$ est grand en module, pr\`es de l'axe r\'eel, et $\chi$ est une fonction r\'eguli\`ere \`a support compact. De tels r\'esultats remontent aux travaux de C. S. Morawetz \cite{Mor1},\cite{Mor2}, qui en d\'eduisait la d\'ecroissance uniforme de l'\'energie locale de l'\'equation d'onde correspondante. 
Il est possible de d\'emontrer (\ref{ineg HF standard}) dans un cadre g\'en\'eral, en modifiant la m\'etrique d\'efinissant le Laplacien ou en rajoutant un obstacle, moyennant une hypoth\`ese essentielle de non-capture sur les g\'eod\'esiques de cette m\'etrique (cf \cite{LaPh},\cite{NB1}, \cite{VZ}).\par
Dans \cite{RV1} et \cite{RV2}, les auteurs consid\`erent des potentiels peu
r\'eguliers, et d\'emontrent des estimations sur la r\'esolvante de $P$, des effets r\'egula\-risants sur les \'equations d'onde et de Schr\"odinger. Le principe de ces deux articles est d'\'ecrire des estimations sur le laplacien libre $\Delta$, puis de consid\'erer $V$ comme une petite perturbation de ce dernier. Un tel raisonnement fonctionne en particulier si $V\in L^{q}$, $q\geq d/2$.\par
Ici, nous supposerons que le potentiel $V$, est petit \`a l'infini et born\'e en dehors d'un ensemble fini de p\^oles distincts:
$$ \PPP=\{p_1,..,p_N\} \subset \RR^d $$ 
pr\`es desquels:
$$ V(x)\approx \frac{a_j}{|x-p_j|^2}. $$
De telles singularit\'es sont critiques, car du m\^eme ordre que le Laplacien. Les singularit\'es moins fortes rentrent dans le cadre de l'article pr\'ecit\'e et pour les singularit\'es d'ordre sup\'erieur on ne peut pas, en g\'en\'eral, d\'emontrer (\ref{ineg HF standard}) (cf \cite{CEx}). Lorsque les constantes $a_j$ sont petites, le potentiel $V$ reste en un certain sens inf\'erieur au laplacien et ce cas est encore trait\'e, lorsque $d\geq 3$, dans \cite{RV2}. En omettant cette hypoth\`ese de petitesse, on change la nature du probl\`eme car on ne peut plus consid\'erer $V$ comme une perturbation du laplacien.\par
 Le cas unipolaire, essentiellement:
$$ P_a=-\Delta+\frac{a}{|x|^2},\quad a+\left(\frac d2-1\right)^2>0,$$
est trait\'e dans \cite{PSTZ} et \cite{BPSTZ}. Moyennant l'hypoth\`ese sur
$a$, qui assure la positivit\'e de l'op\'erateur $P_a$, les auteurs
d\'emontrent des in\'egalit\'es de Strichartz sur les \'equations
d'\'evolutions associ\'ees \`a $P_a$. Leur raisonnement repose de mani\`ere essentielle sur le
caract\`ere radial de $P_a$ et ne se g\'en\'eralise pas au cas
multipolaire.\par
On d\'emontre ici (\ref{ineg HF standard}) pour un potentiel multipolaire. 
Ce type de potentiels appara\^ \i t dans certains mod\`eles de relativit\'e g\'en\'erale, mais la motivation principale de ce travail est l'\'etude d'un probl\`eme critique, cas limite o\`u les singularit\'es de $V$ sont exactement du degr\'e d'homog\'en\'eit\'e du Laplacien.\par
On suppose que pr\`es de chaque p\^ole $p_j$, $V$ est radial (c'est \`a dire fonction de la seule variable $|x-p_j|$). Cette condition est l\'eg\`erement assouplie dans la section \ref{sec.theo2} (cf theor\`eme \ref{theo2}). On fait \'egalement, pr\`es de $p_j$, les hypoth\`eses suivantes:
$$ \frac{a}{|x-p_j|^2}\leq V(x) \leq \frac{C}{|x-p_j|^2},\qquad |\nabla V(x)|\leq \frac{C}{|x-p_j|^3},$$ 
pour une grande constante $C$ et un r\'eel $a$ tel que $a+(d/2-1)^2>0$. Ces hypoth\`eses sont v\'erifi\'es par exemple lorsque $V$ est exactement, pr\`es de chaque $p_j$, de la forme:
$$\frac{a_j}{|x-p_j|^2},\quad a_j+\left( \frac d2 -1 \right )^2 >0. $$
On suppose aussi pour simplifier $V$ nul \`a l'infini et born\'e en dehors des p\^oles. Dans ces conditions on peut toujours d\'efinir, au sens des formes quadratiques, un op\'erateur auto-adjoint semi-born\'e inf\'erieurement $P=-\Delta+V$. Les hypoth\`eses sur $V$ et la construction pr\'ecise de $P$ sont explicit\'ees et discut\'ees dans la section $2$.\par
\begin{theo}
\label{theo1}
    Soit $V$ v\'erifiant les hypoth\`eses (\ref{HypV0}),..,(\ref{HypV5}) et $P=-\Delta+V$. On se donne $\chi \in C^{\infty}_0(\RR^d).$
    Alors:
    \begin{equation}
        \label{ineg1}
        \exists \lambda_{0}>0,\; \exists C>0,\; \forall \lambda >\lambda_{0},\; \forall 
        \eps>0,\;\quad ||\chi R_{\lambda\pm i\eps}\chi||_{L^2\rightarrow 
        L^2}\leq \frac{C}{\sqrt{\lambda}}.
    \end{equation}
Dans le cas d'un seul p\^ole ($N=1$), (\ref{ineg1}) reste vraie sans la borne (\ref{HypV3}) sur la d\'eriv\'ee de $V$.
\end{theo}
Comme d\'ej\`a indiqu\'e, l'estimation (\ref{ineg1}) est l'estimation standard sur la r\'esolvante d'un Laplacien induisant une m\'etrique non captive. 
 Dans les cas captifs, (\ref{ineg1}) est fausse. Le th\'eor\`eme {\ref{theo1} montre en particulier que l'\'energie ne se concentre pas \`a haute fr\'equence sur les p\^oles, et que ces derniers ne se comportent pas non plus comme des obstacles qui renverraient les rayons optiques.\par
On peut d\'eduire directement de (\ref{ineg1}), comme dans \cite{BGT} l'effet r\'egularisant classique sur l'\'equation de Schr\"odinger. L'in\'egalit\'e uniforme haute fr\'equence, telle qu'elle est d\'emontr\'ee ici, implique cet effet localement en temps. Pour avoir des r\'esultats globaux, il faudrait d\'emontrer la m\^eme in\'egalit\'e pour tout r\'eel $\lambda$. On peut aussi d\'eduire du th\'eor\`eme \ref{theo1} la d\'ecroissance de l'\'energie locale de l'\'equation des ondes associ\'ee \`a $P$. Notons aussi que ces r\'esultats impliquent des in\'egalit\'es de Strichartz (avec pertes \'eventuelles) sur les m\^emes \'equations.\par
 La d\'emonstration suit celle de N. Burq \cite{NB1}, et repose sur l'introduction dans un raisonnement par l'absurde d'une mesure de d\'efaut semi-classique, objet introduit ind\'ependamment par P. G\'erard et P.L. Lions (cf \cite{PG}, \cite{LP}). La difficult\'e nouvelle repose dans le compr\'ehension du comportement de la mesure pr\`es de chacun des p\^oles.\par
On peut signaler une approche diff\'erente mais relat\'ee \`a la notre pour d\'emontrer des r\'esultats similaires, celle du calcul de commutateurs positifs introduit par E. Mourre
\cite{Mou}, et qui fonctionne telle quelle sur l'op\'erateur unipolaire
$P_a$. L'auteur tient \`a remercier C. G\'erard et F. Nier pour l'avoir
\'eclair\'e sur ce sujet. On renvoie \`a \cite{CG} (et \cite{ABG} pour
une pr\'esentation g\'en\'erale). L'article \cite{VZ} de A. Vasy et
M. Zworsky  donne une version micro-locale de ce type de techniques.\par
Mentionnons enfin les travaux r\'ecents \cite{VP}, \cite{DAP} sur
l'\'equation des ondes en dimension $3$ avec un potentiel
singulier. L'hypoth\`ese faite dans \cite{DAP} est simplement une
hypoth\`ese de petitesse sur la partie n\'egative du potentiel, similaire
\`a notre hypoth\`ese (\ref{HypV2}). Dans cet article, le potentiel est
pris dans une classe de Kato critique, contenant strictement
$L^{d/2}$. L'introduction de \cite{VP} pr\'esente de mani\`ere tr\`es
compl\`ete les r\'esultats connus sur les op\'erateurs de la forme $-\Delta+V$.\par
La deuxi\`eme partie du texte est consacr\'ee \`a la d\'efinition pr\'ecise de $P$. La troisi\`eme partie concerne la d\'emonstration du th\'eor\`eme \ref{theo1}. Enfin, dans la quatri\`eme partie, on \'enonce et on d\'emontre un raffinnement du th\'eor\`eme \ref{theo1}, o\`u la condition de radialit\'e sur $V$ pr\`es de chaque p\^ole est assouplie. 

\section{D\'efinitions et hypoth\`eses.}\label{sec2}
On commence par expliciter les hypoth\`eses sur $V$. Comme d\'ej\`a pr\'ecis\'e, on se place en dimension $d\geq 2$, et on se donne un ensemble de $N$ p\^oles ($N\geq 1$) distincts:
$$ \PPP=\{ p_1,..,p_N \} \subset \RR^d.$$
On suppose:
\begin{gather}
   \label{HypV0}
    V\text{ \`a support compact sur }\RR^d\\
    \label{HypV1}
    V\in L^\infty_{loc}(\RR^d\backslash \PPP,\RR), 
\end{gather}
et qu'il existe des constantes $l$, $C_V$, $C_V'$ strictement positives, une constante r\'eelle $a$, des fonctions $V_j \in L^{\infty}_{\loc}(]0,l],\RR)$ telles que pour tout entier $j$ compris entre $1$ et $N$, et pour $|x-p_j|\leq l$:
\begin{gather}
\label{HypV2}
     V(x) \geq \frac{a}{|x-p_j|^2},\quad a+\left(\frac{d}{2}-1\right)^2>0\\
\label{HypV3}
\left| V(x) \right| \leq \frac{C_V}{|x-p_j|^2}\\
\label{HypV4}
\left| \nabla V(x) \right| \leq \frac{C_V'}{|x-p_j|^3}\\
\label{HypV5}
V(x)=V_j(|x-p_j|)
\end{gather}
Les hypoth\`eses (\ref{HypV0}) et (\ref{HypV1}) sont loin d'\^etre 
minimales mais on s'int\'eresse ici \`a l'effet des p\^oles sur le comportement de $P$. Pour des potentiels plus g\'en\'eraux \`a l'infini on pourra consulter \cite{BAK}, \cite{VZ}.\par 
La borne (\ref{HypV4}) sur la d\'eriv\'ee et l'hypoth\`ese de radialit\'e (\ref{HypV5}) de $V$ pr\`es des p\^oles sont des hypoth\`eses apparemment techniques. Dans le cas d'un seul p\^ole ($N=1$), (\ref{HypV4}) est inutile, ce qui g\'en\'eralise les r\'esultats connus jusqu'alors pour un potentiel unipolaire radial, qui n\'ecessitait toujours une hypoth\`ese sur la d\'eriv\'ee de $V$. Dans le th\'eor\`eme \ref{theo2} plus loin, on fait une hypoth\`ese l\'eg\`erement plus faible que (\ref{HypV5}), ce qui permet d'inclure des p\^oles de la forme: $1/|x|^2 a\left(x/|x|\right)$. Mais en dehors de ces deux cas, il semble impossible de se passer de (\ref{HypV4}), (\ref{HypV5}) avec la m\'ethode de preuve employ\'ee ici.\par 
Enfin, (\ref{HypV2}) et (\ref{HypV3}) sont absolument essentielles. Sans
l'in\'egalit\'e (\ref{HypV2}), $P$ ne serait plus semi-born\'e inf\'e\-rieurement et le probl\`eme serait d'une nature compl\`etement diff\'erente. La d\'efinition de $P$ comme op\'erateur auto-adjoint, que l'on fait ici par l'extension de Friedrichs, serait elle-m\^eme ambigu\"e. D'autre part, il existe des potentiels unipolaires $V$ strictement positifs, v\'erifiant une majoration juste un peu plus faible que (\ref{HypV3}) (en $\log^2|x|/|x|^2$), et tel que ni (\ref{ineg HF standard}), ni aucune in\'egalit\'e de Strichartz ou de dispersion non triviale ne soient vraie (cf \cite{CEx}). Si l'on veut encore que (\ref{ineg HF standard}) soit v\'erifi\'ee pour des potentiels admettant des p\^oles d'ordre sup\'erieur, il faut donc faire des hypoth\`eses suppl\'ementaires, peut-\^etre supposer une propri\'et\'e monotonie de $V$ au voisinage des p\^oles.\par
Soit $Q$ la forme quadratique sur $L^2$ d\'efinie par:
\begin{equation*}
   \begin{aligned}
       D(Q)=& H^1\cap \left\{u \in L^2(\RR^d); \frac{1}{|x-p_{j}|} u \in 
       L^2(\RR^d), j=1,\dots,N\right\}\\
       Q(u)=& \int |\nabla u|^2 \d x+\int V \abs{u}^2 \d x.
    \end{aligned}   
\end{equation*}
\begin{lem}
\label{lem.Hardy}
Supposons (\ref{HypV0}),..,(\ref{HypV3}). La forme quadratique $Q$ est ferm\'ee, semi-born\'ee inf\'erieurement.
\end{lem}
\begin{proof}
Le fait que $Q$ soit semi-born\'ee inf\'erieurement, trivial pour $d=2$ ($V$ est alors positif pr\`es de chaque p\^ole), d\'ecoule, pour $d\geq 3$, de l'in\'egalit\'e de Hardy:
\begin{equation}
\label{Hardy}
\int |\nabla u|^2 \d x\geq (\frac d2 -1)^2 \int \frac{1}{|x|^2} |u|^2 
\d x,\; u\in H^1(\RR^d).
\end{equation}
En effet, on se donne une partition de l'unit\'e adapt\'ee aux p\^oles $p_j$:
\begin{gather*}
 \chi_0 \in C^{\infty}(\RR^d), \; \chi_j\in C^{\infty}_0(\RR^d),\; \sum_{j=1}^{N} \chi_j^2=1\\
p_k \notin \supp \chi_j,\; k\neq j,\; k=1,\dots,N,\;j=0,\dots,N;
\end{gather*}
et on \'ecrit:
\begin{align*}
Q(u)=&\int \sum_{j=0}^N \chi_j^2 |\nabla u|^2 \d x+\int \sum_{j=0}^{N} V \chi_j^2 |u|^2 \d x\\
=& \int |\chi_0 \nabla u|^2 \d x +\sum_{j=1}^N \int|\nabla(\chi_j u )+[\chi_j,\nabla ]u|^2 \d x + \sum_{j=0}^N \int V|\chi_j u|^2 \d x \\
\geq &
\int |\chi_0 \nabla u|^2 \d x+(1-\eps)\sum_{j=1}^N \int|\nabla(\chi_j u )|^2 \d x\\ 
& -C_{\eps} \sum_{j=1}^N \int|[\chi_j,\nabla]u|^2 \d x + \sum_{j=0}^N \int V|\chi_j u|^2 \d x,
\end{align*}  
o\`u $C_{\eps}$ d\'esigne une constante d\'ependant de $\eps$. On conclut avec l'hypoth\`ese (\ref{HypV2}), et l'in\'egalit\'e (\ref{Hardy}), appliqu\'ee \`a $\chi_j u$, en prenant $\eps$ assez petit. \par
Il est \'evident que $D(Q)$ est complet pour la norme:
$$ ||u||_Q=\sqrt{Q(u)+||u||_{L^2}},$$
 c'est \`a dire que $Q$ est ferm\'ee (en fait, d\`es que $d\geq 3$, $D(Q)$ est exactement l'espace $H^1$ d'apr\`es l'in\'egalit\'e de Hardy). 
\end{proof}
\begin{corol}
Sous les hypoth\`eses (\ref{HypV0}),..,(\ref{HypV3}), on peut associer \`a $Q$ un unique op\'erateur $P$, auto-adjoint, semi-born\'e inf\'erieurement, tel que:
\begin{gather*}
    D(P)=\{u \in D(Q),\; v \mapsto Q(u,v) \text{ continu } L^2\}\\
    \forall u\in D(P),\; \forall v \in D(Q),\; Q(u,v)=(Pu,v)_{L^2}.
\end{gather*}
(cf \cite{RS1})
\end{corol}

La proposition suivante est cons\'equence imm\'ediate de la 
d\'efinition de $P$ et, pour le $1)$, de la densit\'e de 
$C^\infty_{0}(\RR^d\backslash \PPP)$ dans D(Q):
\begin{prop}
\label{D(P)}
    \begin{itemize}
        \item[1)] $D(P)=\{ u\in D(Q); -\Delta u+V u \in L^2\} $,
        \item[2)] Sur $D(P)$,
        $$Pu=-\Delta u+Vu,$$
    \end{itemize}
    o\`u dans $1)$ et $2)$, $-\Delta u+V u$ est \`a prendre au sens 
    des distributions sur $\RR^d \backslash \PPP$.
\end{prop}
\begin{remarque}
M\^eme si $P$ est la somme formelle des op\'erateurs 
$-\Delta$ et $V$, on n'a pas toujours les inclusions:
\begin{gather*}
    D(P)\subset D(-\Delta)=H^2\\
    D(P)\subset D(V)= \{ u\in L^2,\; \frac{1}{|x-p_{j}|^2} u \in 
    L^2,\;j=1,\dots,N\}.
\end{gather*}
Supposons par exemple, que pr\`es d'un p\^ole $p_j$, $V$ soit exactement de la forme:  
\begin{equation}
\label{forme.aj}
\frac{a_j}{|x-p_j|^2},\quad 0<a_{j}+(\frac{d}{2}-1)^2<1. 
\end{equation}
Consid\'erons la fonction:
\begin{equation}
\label{u_s}
u_{s}(x)=|x-p_{j}|^s\varphi(x-p_{j}),\quad 
s=-(d/2-1)+\sqrt{(d/2-1)^2+a},
\end{equation}
o\`u $\varphi$ est \`a support compact dans $\{r<l\}$ et vaut 
$1$ pr\`es de 0. Comme $-\Delta u+Vu$ est nul pr\`es du j-i\`eme 
p\^ole, la proposition pr\'ec\'edente montre que $u$ est un 
\'el\'ement de $D(P)$. Mais $u$ n'est ni dans $H^2$, ni dans 
$D(V)$.
\end{remarque}
\begin{remarque}[non-unicit\'e de l'extension auto-adjointe]
Signalons une autre pathologie int\'eressante, qui survient encore par exemple lorsque $V$ est de la forme (\ref{forme.aj}) pr\`es d'un p\^ole. La fonction $u_{s}$ d\'efinie par (\ref{u_s}), 
avec cette fois:
$$ s=-(d/2-1)-\sqrt{(d/2-1)^2+a},$$
est dans $L^2$ et v\'erifie: 
$$-\Delta u_s+V u_s\in L^2$$
au sens des distribution sur $\RR^d\backslash \PPP$, mais n'est pas 
dans $D(P)$. Si $A$ est l'op\'erateur $-\Delta+V$, d\'efini de 
mani\`ere naturelle sur 
$C_0^{\infty}\left(\RR^d\backslash\PPP\right)$, on a donc:
$$ A\subsetneq P \subsetneq A^*.$$ 
Ainsi l'op\'erateur $A$ n'est pas essentiellement auto-adjoint, et admet 
plusieurs extensions auto-adjointes, dont l'extension de Friedrichs, 
choisie ici, qui est la seule dont le domaine est inclus dans 
$H^1$.
\end{remarque}
Lorsque la constante $a$ de (\ref{HypV2}) est assez grande, les deux remarques pr\'ec\'edentes ne sont plus valables: l'op\'erateur $A$ ci-dessus admet une seule extension auto-adjointe, dont le domaine est exactement l'intersection de $H^2$ et $D(V)$.
Pour un aper\c cu de ses questions et une \'etude du cas unipolaire, 
on pourra consulter \cite[chap X]{RS2}.\par

\section{D\'emonstration du th\'eor\`eme \ref{theo1}}
\label{sec3}
\subsection{Introduction de la mesure et \'etude en dehors des 
p\^oles.}
\label{intro.mu}
On ne fait dans ce paragraphe \ref{intro.mu} que les hypoth\`eses (\ref{HypV0}),..,(\ref{HypV3}) sur le potentiel $V$. L'adjoint de l'op\'erateur born\'e sur $L^2$: $\chi R_{\lambda+i\eps} \chi$ est l'op\'erateur: $\overline{\chi} R_{\lambda-i\eps} \overline{\chi}$. Il suffit donc de d\'emontrer (\ref{ineg1}) avec le signe $-$ devant $i\eps$. Donnons nous une fonction $\chi_1$ de $C^{\infty}_0(\RR^d)$ valant $1$ sur le support de $\chi$. L'in\'egalit\'e (\ref{ineg1}) d\'ecoule de:
\begin{equation}
\tag{\ref{ineg1}'}
\label{ineg2}
\exists \lambda_{0}>0,\; \exists C>0,\; \forall \lambda >\lambda_{0},\; \forall         \eps>0,\;\quad ||\chi_1 R_{\lambda- i\eps}\chi||_{L^2\rightarrow 
        L^2}\leq \frac{C}{\sqrt{\lambda}}.
\end{equation} 
Supposons que (\ref{ineg2}) soit fausse. Alors il existe des suites 
$(\lambda_{n})_{n\geq 1}$, $(\eps_{n})_{n\geq 0}$ telles que:
\begin{gather*}       
\lambda_{n},\eps_n>0,\;\lambda_n\underset{n\rightarrow+\infty}{\longrightarrow}+\infty,\\
    \norme{\chi_{1} R_{\lambda_{n}-i\eps_{n}} \chi }_{\Ld \rightarrow 
\Ld} 
    >\frac{n}{\sqrt{\lambda_{n}}}.
\end{gather*}
Il existe donc une suite $(g_{n})$ de fonctions $\Ld$ telle que:
\begin{equation}
    \label{eq:}
    1=\norme{\chi_{1} R_{\lambda_{n}-i\eps_{n}} \chi g_{n}}_{\Ld} > 
    \frac{n}{\sqrt{\lambda_{n}}} \norme{g_{n}}.
\end{equation}
De plus, on peut supposer $g_{n}\in C^\infty_{0}(\RR^d\backslash 
\PPP)$, cet espace \'etant dense dans $L^2$. 
On se place dans un contexte semi-classique, en posant:
$$ u_{n}=R_{\lambda_{n}-i\eps_{n}} \chi g_{n},\; 
h_{n}=\frac{1}{\sqrt{\lambda_{n}}},\; f_{n}=h_{n} g_{n},\; 
\alpha_{n}=\eps_{n} h_{n}.$$
On a donc:
\begin{gather}
    \label{eq sur u}
  h_n^2(-\Delta+V) u_{n}-(1-i\alpha_{n}h_{n})u_{n}=\chi h_{n} f_{n},\\
\label{hyp u}
        \norme{\chi_{1} u_{n}}_{\Ld}=1,\; \norme{f_{n}}_{\Ld}
        \underset{n\rightarrow +\infty}{\longrightarrow} 0,\; 
        \alpha_{n}>0. 
\end{gather}
Notons que contrairement \`a un ``vrai'' probl\`eme semi-classique, 
le 
potentiel $V$ a un coefficient $h_{n}^2$. Ce coefficient le rend 
n\'egligeable, sauf au voisinage de chaque p\^ole.\par
On omettra parfois les indices $n$ pour all\'eger les 
notations. On commence par noter que l'on peut se contenter d'\'etudier la r\'esolvante au voisinage de l'axe r\'eel:
\begin{lem}
    \begin{gather}
        \label{alpha 0}
        \alpha_n \underset{n\rightarrow +\infty}{\longrightarrow} 0,\\
        \label{alpha 0 bis}
        \alpha_n ||u_n||^2_{\Ld} \underset{n\rightarrow +\infty}{\longrightarrow} 0.
    \end{gather}
\end{lem}
\begin{proof}
Ceci d\'ecoule imm\'ediatement du caract\`ere auto-adjoint de 
l'op\'erateur $P$. En effet, en multipliant (\ref{eq sur u}) par 
$\ubar_n$ et en 
int\'egrant la partie imaginaire, on obtient:
$$ \Im \big(h^2 Q(u_n)-\|u_n\|^2_{L^2}\big)+h_n\alpha_n 
\|u_n\|^2_{L^2}=h\Im\int \chi_{2} f_n \ubar_n\,dx $$
On en d\'eduit (\ref{alpha 0 bis}) car d'apr\`es (\ref{hyp u}), 
$\chi u_n$ est born\'e et 
$f_n$ tend vers $0$ dans $L^2$. Egalement d'apr\`es 
(\ref{hyp u}), la norme de $u_n$ dans 
$L^2$ est minor\'ee quand $n$ tend vers $+\infty$ par un r\'eel strictement 
positif, ce qui donne (\ref{alpha 0}).
\end{proof}
Le lemme suivant, cons\'equence facile de l'in\'egalit\'e de r\'esolvante (\ref{ineg1}) sur le Laplacien libre, est 
d\'emontr\'e dans \cite{NB1}:
\begin{lem}
La suite $(u_{n})$ est born\'ee dans $L^2_{loc}(\RR^d)$.
\end{lem}
On peut donc introduire la mesure semi-classique $\mu$ 
associ\'ee \`a $u_{n}$ et \`a la suite d'\'echelles $(h_{n})$ (cf 
\cite{NB1}, \cite{CFK}). C'est 
une mesure de Radon positive sur $\RR^d_{x}\times\RR^d_{\xi}$ qui 
v\'erifie, \`a extraction d'une sous-suite de $(u_n,h_n)$ pr\`es:
\begin{equation}
\label{lim def mu}
 \forall a \in \Cio(\RR^d\times\RR^d),\quad \big(a(x,h_nD)\varphi(x) 
u_n,u_n\big) 
\underset{n\rightarrow +\infty}{\longrightarrow} <\mu,a>
\end{equation}
O\`u l'on a not\'e $\big(.,.\big)$ le produit scalaire $\Ld$, 
$\varphi\in C^\infty_0(\RR^d_x)$ une fonction valant $1$ sur la 
projection en $x$ du support de $a$, et $a(x,h_n D)$ la suite d'op\'erateur 
uniform\'ement born\'ee sur $\Ld$ de 
noyau:
$$  \surpid \int a(x,h_n\xi) e^{i(x-y).\xi}\d\xi. $$
Rappelons que la limite (\ref{lim def mu}) ne depend pas de la 
fonction $\varphi$.
\begin{prop}
\label{prop.mu}
Sous les hypoth\`eses (\ref{HypV0}),..,(\ref{HypV3}):
\begin{enumerate}
\item \label{h-oscillante} {\bf Vitesse d'oscillation de $u_n$ et convergence $L^2$.} Si $\psi \in C_0^{\infty}(\RR^d)$,
\begin{equation}
\label{nabla.u}
h_n^2 \int |\nabla u_n |^2 \psi \d x+\sum_{j=1..N} h_n^2\int \frac{1}{|x-p_j|^2} |u_n|^2 \psi \d x=O(1),\quad n\rightarrow+\infty.
\end{equation}
En particulier:
\begin{equation}
\label{equivalent.mu}
u_n \underset{n\rightarrow +\infty}{\longrightarrow} 0\text{ dans } L^2_{\loc}(\RR^d) \iff \mu=0.
\end{equation}
\item {\bf Localisation de $\mu$}.
\label{localisation.mu}
    Le support de la mesure $\mu \indic_{(\RR^d\backslash \PPP)\times 
    \RR^d}$ est inclus dans $\{(x,\xi)\in 
\RR_{x}^d\times\RR_{\xi}^d,\; 
    \carre{\xi}=1\}$.
\item {\bf Invariance de $\mu$}.
\label{invariance.mu}
    La mesure $\mu$ v\'erifie l'\'equation:
\begin{equation}
\label{transport}    
\xi.\partial_{x}\mu=0,
\end{equation} 
    au sens $\DDD'((\RR_{x}^d\backslash \PPP)\times 
    \RR_{\xi}^d)$.
\item {\bf Condition \`a l'infini.} \label{incoming}
  Soit $M>0$ assez grand.
    La mesure $\mu$ est nulle pr\`es des points rentrants:
    $${\rm Inc}=\{(x,\xi)\in \RR^d\times \RR^d;\; \abs{x}\geq M,\;x.\xi\leq 0\}. $$
\end{enumerate}
\end{prop}
\begin{proof}
Pour d\'emontrer le \ref{h-oscillante}, on fait le produit scalaire de l'\'equation (\ref{eq sur u}) sur $u$ avec $\psi u$ et on obtient, en utilisant que $f_n$ et $u_n$ sont born\'es dans $L^2_{\loc}$:
$$ h_n^2 \Re Q(u_n,\psi u_n)=O(1),\quad n\rightarrow +\infty.$$
Mais par une int\'egration par partie \'el\'ementaire:
$$ \Re Q(u_n,\psi u_n)=h_n^2 \int |\nabla u_n|^2\psi \d x-h_n^2 \int |u_n|^2 \Delta\psi \d x+ h_n^2 \int V |u_n|^2 \psi \d x,$$
ce qui implique (\ref{nabla.u}), en utilisant l'in\'egalit\'e de Hardy comme dans le lemme \ref{lem.Hardy}. On en d\'eduit que $u_n$ est $h_n$-oscillante, c'est \`a dire que pour toute fonction $\psi$ de $C^{\infty}_0(\RR^d)$,
\begin{equation}
\label{def.h-oscillante}
\lim_{R\rightarrow  +\infty} \limsup_{n\rightarrow +\infty} \int_{|h_n \xi|\leq 1} |\widehat{\psi u_n}(\xi)|^2 \d \xi=0.
\end{equation}  
L'\'equivalence (\ref{equivalent.mu}) est une cons\'equence facile de
(\ref{def.h-oscillante}) (cf \cite{CFK})
\par
Les points \ref{localisation.mu} et \ref{invariance.mu} sont
\'el\'ementaires et d\'emontr\'es dans \cite{NB1}
Le point \ref{incoming}, moins imm\'ediat, se d\'eduit de l'\'etude du
Laplacien libre. C'est une version micro-locale de la condition de
radiation Sommerfeld
. Si l'on change le signe devant $i h \alpha_n$ dans (\ref{eq sur u}), il faut remplacer ${\rm Inc}$ par 
   $${\rm Outc}=\{(x,\xi);\; \abs{x}\geq M,\;x.\xi\geq 0\}. $$
On renvoie \`a \cite{NB1}
pour la d\'emonstration.    
\end{proof}
\newlength{\debutpar}
\setlength{\debutpar}{\parindent}

\begin{minipage}{60mm}
\dessmua
\begin{center}
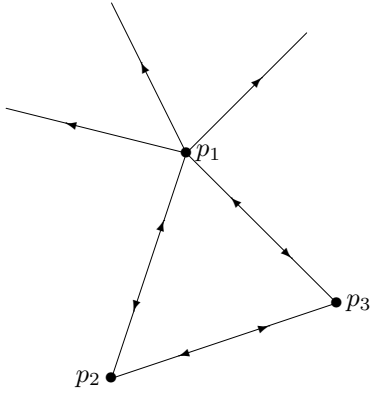
Fig. 1: le support de $\mu$
\end{center}
\end{minipage}
\begin{minipage}{95mm}
\setlength{\parindent}{\debutpar}
L'\'equation (\ref{transport}) est une \'equation de transport qui implique que la mesure, en dehors des p\^oles, est invariante le long des courbes int\'egrales du champ hamiltonien associ\'e \`a $\xi \partial_{x}$, qui dans notre cas sont de la forme: $ \{ (x_0+s\xi_0,\xi_0), s\in ]a,b[ \}$. Dans \cite{NB1}, il n'y a pas de p\^ole, et le point \ref{incoming} implique donc, avec cette propri\'et\'e d'invariance de la mesure et une hypoth\`ese de non-capture qui dit que toute courbe int\'egrale du champ hamiltonien passe dans l'ensemble ${\rm Inc}$, que la mesure $\mu$ est nulle. Ceci montre, la suite $(u_n)$ \'etant $h_n$-oscillante, qu'elle tend vers $0$ dans $L^2_{\loc}$, contredisant ainsi l'hypoth\`ese (\ref{hyp u}).\par

Dans notre cas, la strat\'egie de preuve est la m\^eme, mais l'argument pr\'ec\'edent ne fonctionne pas pour les trajectoires passant par $\PPP$. L'invariance de la mesure et le point \ref{incoming} de la proposition \ref{prop.mu} impliquent seulement que le support de $\mu$ est inclus dans la r\'eunion des rayons reliant les p\^oles et des rayons sortants partant de chacun de ces p\^oles (cf figure ci-jointe).
\end{minipage}

Ces derniers rayons sont les plus faciles \`a \'eliminer, par un argument de ``conservation de l'\'energie'' exprim\'e dans la proposition \ref{mu a l'infini}: si le support de la mesure $\mu$  ne contient aucun rayon rentrant dans un certain compact, il ne peut pas non plus contenir de rayon sortant de ce compact.

\begin{prop}
    \label{mu a l'infini}
    Soient $R_{1},R_{2}$, v\'erifiant: 
    $$0<R_{1}<R_{2},$$ 
    et tels que $\chi$ et $V$ soient nuls sur $\{ |x| \geq 
R_{1}\}$. Supposons: 
    \begin{equation}
        \label{x.xi>0}
        \text{supp}\,\mu\cap \{ R_{1}\leq |x|\leq R_{2} \} \subset 
        \{x.\xi > 0\}
    \end{equation}
    ou:
    \begin{equation}
        \tag{\ref{x.xi>0}'}
        \label{x.xi<0}
        \text{supp}\,\mu\cap \{ R_{1}\leq |x|\leq R_{2} \} \subset 
        \{x.\xi < 0\}
    \end{equation}
    Alors $\mu$ est nulle sur $\{|x|\geq R_{1}\}$.
\end{prop}
\begin{proof}
Soit $\varphi\in C^\infty(\RR^d)$, radiale, telle que:
$$ r\geq R_{2} \Rightarrow \varphi(r)=1,\quad r\leq R_1 \Rightarrow 
\varphi(r)=0, \quad \frac{\partial\varphi}{\partial r}\geq 0. $$
Alors:
\begin{align*}
    0&=\Im \big( (-h^2\Delta u-u+ih\alpha u),\varphi u\big)_{\Ld}\\
    &=\Im h \int h\nabla \varphi.\nabla u \ubar\,dx+h\alpha \int 
|u|^2\varphi \d x\\
\end{align*}
Donc d'apr\`es (\ref{alpha 0 bis}):
\begin{equation}
\label{grad tend vers 0}
\Im \int h_n \nabla \varphi.\nabla u_n \ubar_n \d x \underset{n \rightarrow 
+\infty}{ \longrightarrow} 0. 
\end{equation}
La suite $u_n$ est $h_n$-oscillante, on a donc, en se donnant une fonction $\psi$ de $C_0^{\infty}(\RR^d)$, \`a valeurs r\'eelles et valant $1$ pr\`es de $0$,
\begin{align*}
 \lim_{n \rightarrow +\infty} \Im \int h\nabla \varphi.\nabla u \ubar \d x
=&\lim_{R \rightarrow + \infty} \lim_{n \rightarrow +\infty} \Im \int 
 h\nabla \varphi.\nabla u\, \psi\left(R^{-1} hD\right)\ubar \d x\\
=&\lim_{R \rightarrow +\infty} \lim_{n\rightarrow +\infty} \Im \int \psi\left(R^{-1} h D\right)\left(\nabla\varphi.\nabla u\right) \ubar \d x\\
=& -\big< \mu,\nabla \varphi(x).\xi \big>. 
\end{align*}
On d\'eduit de (\ref{grad tend vers 0}) que cette derni\`ere 
quantit\'e est nulle.
Mais d'apr\`es l'hypoth\`ese (\ref{x.xi>0}) ou (\ref{x.xi<0}), 
$\nabla \varphi(x).\xi$ a un signe constant sur le support de $\mu$. 
La 
positivit\'e de $\mu$ montre alors qu'elle est nulle sur l'ensemble:
$$ \{ \nabla\varphi(x).\xi\neq 0\}$$
ce qui implique la nullit\'e de $\mu$ sur la couronne $\{\frac{d \varphi}{d r}\neq 0\}$, et donc par 
invariance sur $\{|x|\geq R_{1} \}$.
\end{proof}
Il nous reste \`a montrer que $\mu$ ne charge pas les p\^oles et s'annule le long des rayons reliant ces p\^oles. La premier point ne n\'ecessite pas l'hypoth\`ese (\ref{HypV4}) et est trait\'e dans les deux prochaines parties. Dans le cas d'un seul p\^ole, on en d\'eduit imm\'ediatement la contradiction recherch\'ee. Le deuxi\`eme point est trait\'e dans la partie \ref{multipolaire} o\`u on a besoin de toutes les hypoth\`eses (\ref{HypV0}),.., (\ref{HypV5}) sur $V$.

\subsection{Elimination des petites harmoniques sph\'eriques pr\`es d'un p\^ole.}
\label{petites}
D'apr\`es le paragraphe \ref{intro.mu}, la projection spatiale du support de la mesure $\mu$ est incluse dans l'ensemble form\'e des $N$ p\^oles et des $\frac{N(N-1)}{2}$ segments les reliant. 
On cherche ici \`a \'etudier le comportement de $\mu$ au voisinage d'un p\^ole $p_j$. On translate le rep\`ere pour prendre $p_j$ comme origine et on se place en coordoon\'ees sph\'eriques:
$$ r=|x|\in ]0,+\infty[,\quad \theta=\frac{x}{|x|}\in S^{d-1}.$$ 
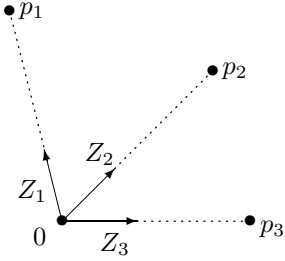
\begin{minipage}{50mm}
\dessmub
\begin{center}
Fig. 2: les vecteurs $Z_j$
\end{center}
\vspace{3mm}
\end{minipage}
\begin{minipage}{105mm}
Notons $\PPP'_0$ l'ensemble des p\^oles diff\'erents de $0$ dans ce nouveau rep\`ere et $Z_j$ les directions sortantes: 
$$Z_{j}=\frac{p_{j}}{|p_{j}|},\quad p_j\in \PPP'_0.$$  
\end{minipage}

La mesure au voisinage de $x=0$ se concentre (sauf peut-\^etre en $x=0$), sur $\{|\xi|=1\}$, sur chacun des segments partant de $0$ dans les directions $Z_{j}$. Elle est invariante par le flot hamiltonien en dehors de $0$ et ne charge pas $0$. On en d\'eduit:
\begin{equation}
    \label{exp de mu}
\begin{gathered}
    \big<\indic_{x\neq 0} \mu,a\big>= \sum_{p_{j}\in\PPP'} \int_{0}^{+\infty} 
    \left(\lambda_{i}^{+} a(x=tZ_{i},\xi=Z_{i})+\lambda_{i}^{-} 
    a(x=tZ_{i},\xi=-Z_{i})\right)\,dt,\\
a\in C_0^{\infty}\left(\{x\in\RR^d,|x|\leq l\}\times\RR^d_{\xi}\right),
\end{gathered}
\end{equation}
o\`u les $\lambda_{i}^{+}$, $\lambda_{i}^{-}$ sont, du fait de la 
positivit\'e 
de $\mu$, des constantes positives. \par
D\'ecomposons $u$ et $f$ en harmoniques sph\'eriques:
\begin{equation}
    \label{harmo}
    u_{n}(x)=\sum_{k\in \NN} \ukn(r)e_k(\theta),\quad 
    f_{n}(x)=\sum_{k\in\NN} \fkn(r) e_k(\theta),
\end{equation}
o\`u les $e_k$ sont les fonctions propres du Laplacien sur la 
sph\`ere, qui forment une base hilbertienne de $L^2(S^{d-1})$ telle 
que:
\begin{gather*}
    e_{k}\in C^{\infty}(S^{d-1}),\quad-\Delta e_k=\nu_{k}^2 e_{k}\\
    \nu_{k+1}\geq \nu_{k}\geq \nu_{0}=0,\quad 
    \nu_{k}\underset{k\rightarrow +\infty}{\longrightarrow}+\infty  
\end{gather*}
Fixons un entier naturel $\tilde{\nu}$. On s\'epare les petites et les 
grandes harmoniques sph\'eriques de $u$ et de $f$:
\begin{equation}
\label{p+g}
\begin{gathered}
    u_{n}=u_{\pp n}+u_{\gg n}, \quad u_{\pp n}\overset{def}{=} 
\sum_{\nu_{k}\leq 
    \tilde{\nu}} \ukn e_k,\quad 
     u_{\gg n}\overset{def}{=} \sum_{\nu_{k}>
    \tilde{\nu}} \ukn e_k\\
    f_n=f_{\pp n}+f_{\gg n}, \quad f_{\pp n}\overset{def}{=} 
\sum_{\nu_{k}\leq 
    \tilde{\nu}} \fkn e_k,\quad 
     f_{\gg n}\overset{def}{=} \sum_{\nu_{k}>
    \tilde{\nu}} \fkn e_k.
\end{gathered}
\end{equation}
On notera $\mu_{\pp}$ (respectivement $\mu_{\gg}$) la mesure semi-classique
associ\'ee pr\`es de $0$ \`a la suite $(u_{\pp n})_n$ (respectivement
$(u_{\gg n})_n$) et \`a la suite d'\'echelle $(h_n)_n$. Par un argument
\'el\'ementaire d'orthogonalit\'e, l'\'equation (\ref{eq sur u}) \'etant
radiale au voisinage du p\^ole, les deux suites v\'erifient pr\`es de
$0$, cette m\^eme \'equation (en rempla\c cant $f_n$ par $f_{\gg n}$
ou $f_{\pp n}$), et elles sont donc toutes les deux $h_n$-oscillantes.\par
L'int\'er\^et de cette d\'ecomposition est que pour \'etudier $u_{\pp}$, on est ramen\'e \`a un nombre fini d'\'equations diff\'erentielles ordinaires, et que l'op\'erateur $P$ est ``tr\`es'' positif pr\`es de $0$ lorsqu'il agit sur les grandes harmoniques sph\'eriques.
Dans ce paragraphe, on montre que les petites harmoniques sph\'eriques ne jouent aucun r\^ole. L'\'etude de $\mu_{\gg}$ se fera dans les deux paragraphes suivants.
\begin{prop}
\label{mup.nulle}
Sous les hypoth\`eses (\ref{HypV0}),..,(\ref{HypV3}) et (\ref{HypV5}), la mesure $\mu_{\pp}$ est nulle et les mesures $\mu$ et $\mu_{\gg}$ sont \'egales.
\end{prop} 
\begin{proof}
On commence par d\'emontrer:
\begin{lem}
\label{lem.mup.0}
La mesure $\mu_{\pp}$ ne charge pas $0$.
\end{lem}

\begin{proof}
On consid\`ere la mesure $\mu_k$ associ\'ee \`a la suite $h_n$-oscillante
$\left(\ukn e_k\right)_n$. Puisque $u_{\pp}$ est la somme d'un nombre fini de
$\uk e_k$, il suffit de montrer le lemme sur chacune des mesure $\mu_k$. On
fixe donc $k\geq 0$.
On a, en notant $'$ la d\'eriv\'ee par rapport \`a $r$:
\begin{equation}
    \label{eq.uk}
    -h_n^2\ukn'' + h_n^2 \frac{d-1}{r} 
    \ukn'  + h_n^2 (V+\frac{\nu_k^2}{2}) \ukn 
    -(1-ih_n\alpha_n) \ukn= h \fkn
\end{equation}
Posons: 
\begin{gather}
\ukn(r)= r^{\frac{d-1}{2}} \vkn,\quad \fkn(r)=r^{\frac{d-1}{2}}
\gkn\\
\label{bornes W_k}
W_k\overset{def}{=}V+\frac{\nu_k^2}{r^2}+\frac{d^2-4d+3}{4r^2},\quad |W_k(r)|\leq \frac{C_1}{r^2}
\end{gather}
On a:
\begin{gather}
\label{eq vk}
S_k \vk=h\gk,\quad S_k\overset{def}{=} -h^2 \frac{d^2}{dr^2} +h^2 W_k-(1-ih\alpha)
\end{gather}
Dans cette d\'emonstration, on ne pr\'ecise pas la d\'ependance
\'eventuelle en $k$ des constantes (qui sont bien s\^ur
ind\'ependantes de $n$), $k$ \'etant fix\'e de bout en bout.
Il suffit de montrer:
$$ \forall \eps>0,\;\exists r_1>0,\; \limsup_{n\rightarrow +\infty} \int_0^{r_1} \left|\vkn\right|^2\,dr \leq \eps.$$ 
Soit $\eps$ strictement positif. D'apr\`es (\ref{nabla.u}),
$$ \exists C_2>0, \quad \forall n,\quad \int_0^{r_0} \frac{h_n^2}{r^2}|\vkn|^2\,dr\leq C_2.$$ 
Soit $m>0$. On a:
\begin{equation}
\label{majoration.E.1}
\int_0^{mh} |\vk(r)|^2\d r\leq m^2 \int_0^{mh} \frac{h^2}{r^2} |\vk|^2\d r\leq C_2 m^2\leq \eps,
\end{equation}
en choisissant $m$ assez petit pour que la derni\`ere in\'egalit\'e soit v\'erifi\'ee. Il nous reste \`a majorer la norme $L^2$ de $\vk$ dans la zone $\{r>mh\}$. On introduit:
$$ E_{k n}(r)\overset{def}{=}\left|\vkn(r)\right|^2+\left|h_n \vkn'(r) \right|^2,$$
qui est d\'erivable, de d\'eriv\'ee:
\begin{align}
\notag
E_k'(r)=&2\Re\left( \vk' \vkbar+h^2 \vk'\vkbar'' \right)\\
\notag
=& 2\Re\left(h^2 W_k  \vk' \vkbar-i \alpha \vk'\vkbar-h\chi \vk' \fkbar\right)\\
\notag
\left|E_k'(r) \right|\leq& \left(\frac{C_1h}{r^2}+\alpha\right)|h\vk'\vk|+|h\vk' \fk|\\
\label{majoration.E'}
-E_k'(r)\leq & \left(\frac{C_1h}{r^2}+1\right)E_k(r)+\left| \fk(r)\right|^2,
\end{align}
d\`es que $n$ est assez grand pour que $\alpha$ soit inf\'erieur \`a $1$. On a obtenu la deuxi\`eme ligne par l'\'equation
(\ref{eq vk}) sur $\vk$. Majorons $E_k$ pr\`es de $0$ par le lemme de
Gronwall. Pour cela, on fixe deux r\'eels strictement positifs $t$ et $\rho$. Il d\'ecoule de (\ref{majoration.E'}):
\begin{equation*}
\forall r\geq t,\;
-\frac{d}{dr}\left( e^{\int_{t}^r (C_1 h/s^2+1)\d s} E_k(r)\right)\leq e^{\int_{t}^r (C_1 h/s^2+1)\d s}|\fk(r)|^2.
\end{equation*}
Soit, en int\'egrant cette in\'egalit\'e entre $t$ et $t+\rho$, avec $t\geq mh$:
\begin{equation*}
E_k(t)\leq e^{C_1/m+\rho} \kappa_k^2+ e^{C_1/m} E_k(t+\rho),\quad \kappa_k^2(n)\overset{def}{=}\int_0^l \left|\fkn(r)\right|^2\d r.
\end{equation*}
On int\`egre maintenant par rapport \`a la mesure $\d t$ entre $mh$ et un r\'eel strictement positif $r_1$. On obtient:
\begin{equation}
\label{majoration.E.2}
\int_{mh}^{r_1} E_k(t) \d t\leq \underbrace{e^{C_1/m+\rho} r_1\kappa_k^2}_{\underset{n\rightarrow +\infty}{\longrightarrow 0}} +e^{C_1/m+\rho} \int_{mh}^{r_1} E_k(t+\rho) \d t
\end{equation}
Fixons $\rho$ strictement compris entre $0$ et $l$. On a:
\begin{gather}
\notag
\int_{\rho}^{r_1+\rho} E_k(r) \d r=\int_{\rho}^{r_1+\rho} |u_k(r)|^2 r^{d-1} \d r+ h^2 \int_{\rho}^{r_1+\rho} |u_k'(r)|^2 r^{d-1} \d r +O(h^2), \quad n\rightarrow +\infty\\
\label{majoration.E.3}
\limsup_{n\rightarrow +\infty}\int_{\rho}^{r_1+\rho}E_{kn}(r) \d r \leq  2\mu_k(\{\rho\leq |x|\leq r_1+\rho\})
\end{gather}  
On a obtenu cette derni\`ere ligne en utilisant que sur le support de $\mu\indic_{\{x\neq 0\}}$, et donc sur celui de $\mu_k\indic_{\{x\neq 0\}}$, $|\xi|$ vaut $1$.
La forme de la mesure $\mu_p$ pr\`es de $0$ montre que le terme de droite de
cet in\'egalit\'e peut-\^etre choisi plus petit que $e^{-C_1/m} \eps$, en prenant
$r_1$ assez petit. En effet, l'orthogonalit\'e des harmoniques sph\'eriques montre que:
\begin{align*}
\mu(\{ \rho\leq |x| \leq r_1+\rho \})=&\mu_{\pp}(\{ \rho\leq |x| \leq r_1+\rho \})+\mu_{\gg}(\{ \rho\leq |x| \leq r_1+\rho \})\\
\geq &\mu_{\pp}(\{ \rho\leq |x| \leq r_1+\rho \}),
\end{align*}
et la mesure $\mu$ ne charge pas les cercles autour du p\^ole.
Finalement en utilisant (\ref{majoration.E.3}) avec
un tel $r_1$, ainsi que
(\ref{majoration.E.1}),(\ref{majoration.E.2}), on obtient: 
\begin{equation*}
\limsup_{n\rightarrow+\infty} \int_0^{r_1} E_{kn}(r)\d r\leq 2\eps,
\end{equation*}
ce qui ach\`eve la d\'emonstration.
\end{proof}

Fixons maintenant un $r_0$ proche de $0$. Les suites $(u_n)$, $(u_{\gg n})$ et $(u_{\pp n})$, v\'erifient deux 
propri\'et\'es d'``orthogonalit\'es''.\par
\begin{itemize}
\item
La premi\`ere est simplement l'orthogonalit\'e dans $L^2$ 
des harmoniques 
    sph\'eriques, qui implique:
    $$ \big( u_{\gg n},u_{\pp n}\big)_{L^2(\{r\leq r_0\})}=0.$$
    D'o\`u:
    \begin{gather}
        \notag
        |u_n|^2_{L^2(\{r\leq r_0\}}=|u_{\pp n}|^2_{L^2(\{r\leq 
        r_0\}}+|u_{\gg n}|^2_{L^2(\{r\leq r_0\}},\\
        \label{ortho 1}
        \mu(\{r\leq r_0\})=\mu_{\gg}(\{r\leq r_0\})+\mu_{\pp}(\{r\leq r_0\}).
    \end{gather}
\item La deuxi\`eme est le fait que les deux mesure $\mu\indic_{\{|x|<2 r_0\}}$ et $\mu_{\pp}\indic_{\{|x|< 2r_0 \}}$ sont mutuellement singuli\`eres. La premi\`ere de ces mesures est port\'ee par l'ensemble:
$$\SSS=S\times \RR^d_{\xi}, \quad S=\left\{ s Z_j,\quad s\in [0,2 r_0[,\quad p_j\in \PPP'_0 \right\}.$$ 
Nous allons d\'emontrer que $\mu_{\pp}\indic_{\{|x|< 2r_0 \}}$ ne charge pas cet ensemble. Soit $\eps>0$. La mesure $\mu_p$ ne charge pas $0$, donc si $\rho$ est assez petit,
\begin{equation}
\label{charge.pas}
\mu_{\pp}(\{|x|\leq \rho\}) \leq \eps.
\end{equation}
D'autre part, $u_{\pp}$ ne peut pas se concentrer sur un rayon. En effet, si $Z$ est dans $S^{d-1}$:
\begin{equation}
\label{charge.pas2}
\begin{aligned}
\left\{\int_{\big|x/|x|-Z\big|\leq \eta,\; |x|\leq 2r_0} \big|u_{\pp n}(x)\big|^2 dx\right\}^{1/2}\leq& \sum_{k,\nu_k\leq \tilde{\nu}} \left\{\int_{\big|x/|x|-Z\big|\leq \eta,\;|x|\leq 2r_0} \big|\ukn(|x|) e_k(x/|x|)\big|^2 dx\right\}^{1/2}\\
\leq & \sum_{k} \left\{\int_0^{2r_0} |\ukn(r)|^2 r^{d-1} \d r\right\}^{1/2} \left\{\int_{|\theta-Z|\leq \eta} |e_k(\theta)|^2 \d \theta\right\}^{1/2}
\end{aligned}
\end{equation}
Le premier facteur est major\'e par une constante ind\'ependante de $n$, et le deuxi\`eme, ind\'ependant de $n$, tend vers $0$ lorsque $\eta$ tend vers $0$. On obtient finalement, avec (\ref{charge.pas}) et (\ref{charge.pas2}),  que pour tout $\eps$, il existe un voisinage $V$ de $\SSS$ tel que:
$$ \mu_{\pp}(V) \leq \eps,$$
et donc que $\mu_{\pp}\indic_{\SSS}=0$. Il est facile de montrer, dans ces 
    conditions que pour tout symbole $a(x,\xi)\in 
    C^\infty_{0}(\{|x| <2r_0\}\times \RR^d)$ et pour toute fonction r\'eguli\`ere $\varphi$ tronquant autour de $0$, on a:
\begin{equation}
\label{up.u.0}
    (a(x,hD) \varphi u_{\pp},u)_{L^2} \underset{n\rightarrow 
+\infty}{\longrightarrow} 
    0. 
\end{equation}
(il suffit de diviser le produit scalaire entre un petit voisinage de $S$, sur lequel la limite sup\'erieur des normes $L^2$ de $u_{\pp n}$ est aussi petite que l'on veut, et son compl\'ementaire, sur lequel la norme $L^2$ de $u_{\gg n}$ tend vers $0$ quand $n$ tend vers l'$\infty$).
    Finalement, l'\'egalit\'e: $u_{\gg}=u-u_{\pp}$ et (\ref{up.u.0}) impliquent:
    \begin{equation}
        \label{ortho 2}
        \mu_{\gg}(\{r\leq r_0\})=\mu(\{r\leq r_0\})+\mu_{\pp}(\{r\leq r_0\}).
    \end{equation}
\end{itemize}
D'apr\`es (\ref{ortho 1}) et (\ref{ortho 2}), $\mu_{\pp}(\{r\leq r_0\})$ 
est nulle. La proposition \ref{mup.nulle} est d\'emontr\'ee.
\end{proof}
\subsection{Absence de concentration sur le p\^ole}
\label{concentration}
\begin{prop}
\label{pas.de.concentration}
Sous les hypoth\`eses (\ref{HypV0}),..,(\ref{HypV3}) et (\ref{HypV5}), la mesure
$\mu$ ne charge pas les p\^oles.
\end{prop}
\begin{corol}
\label{mu.nulle.N=1}
Sous les m\^emes hypoth\`eses et si $N=1$, $\mu$ est nulle.
\end{corol}
En effet, dans ce dernier cas, d'apr\`es le paragraphe \ref{intro.mu}
(propositions \ref{prop.mu} et \ref{mu a l'infini}), la projection en $x$ du
support de $\mu$ est incluse dans le p\^ole. L'\'etude du cas $N\geq 2$ est
compl\'et\'ee dans la partie \ref{multipolaire}.\par
Pour montrer la proposition, on reprend les notations de la partie
\ref{petites}. On se place encore au voisinage du p\^ole $p_{j_0}=0$ et on
consid\`ere la d\'ecomposition (\ref{p+g}) en petites et grandes harmoniques
sph\'eriques: $u=u_{\pp}+u_{\gg}$.  La mesure $\mu_{\pp}$
\'etant nulle, il suffit de montrer que la mesure $\mu_{\gg}$ ne charge pas
les p\^oles. On va en fait montrer un r\'esultat plus fort:
\begin{lem}
\label{le.lem.ug.en.0}
Soit $t\in ]0,1[$. Si $\tilde{\nu}$ est choisi assez grand:
\begin{equation}
\label{lem.ug.en.0} \int_{|x|\leq 1} |x|^{-t} |u_{\gg n}|^2 \d x =O(1),\quad
n\rightarrow +\infty. 
\end{equation}
\end{lem}
\begin{proof}[D\'emonstration]
Comme dans la d\'emonstration du lemme \ref{lem.mup.0}, on note
$r^{(d-1)/2}\vkn=\ukn$.  
Fixons $n$ et $k$. Pour justifier nos int\'egrations par parties, on a besoin de conna\^itre le comportement pr\`es de $0$ de $\vkn$ et de ces d\'eriv\'ees. On notera:
$$ F(r) \underset{r \rightarrow 0}{\lesssim} G(r) \iff \exists
\eps>0,\quad\exists A>0,\quad \forall r\in ]0,\eps[,\quad |F(r)| \leq A |G(r)|. $$
On a choisi les fonctions $f_n$ nulles pr\`es de chaque p\^ole. La fonction
$\vkn$ est donc, au voisinage de $0$, solution de l'\'equation
diff\'erentielle en $y$: $S_ky=0$ (o\`u $S_k$ est l'op\'erateur
diff\'erentiel de degr\'e $2$ d\'efini par (\ref{eq vk})), qui admet une
famille libre de solutions $\{y_+,y_-\}$ telle que:
\begin{equation}
\label{equivalences}
\begin{aligned}
    \left(\frac{d}{dr}\right)^j y_{+}(r) &\underset{r\rightarrow 0}{\lesssim} 
    r^{1-j+\sigma_k},\quad j=0,1,2\\
 \left(\frac{d}{dr}\right)^j y_{-}(r) &\underset{r\rightarrow 0}{\gtrsim}  
    r^{1-\sigma_k}\\
 \sigma_k &\egaldef \sqrt{ (d/2-1)^2 +a + \nu_{k}^2 }.
\end{aligned}
\end{equation}
R\'eservons la preuve de cette affirmation \`a plus tard (cf le lemme \ref{lem.equivalence}).
La fonction $u_n$ \'etant dans le domaine de $P$, $r^{-1}\ukn$ est dans
l'espace $L^2(r^{d-1}dr)$, et donc $r^{-1} \vkn$ est dans l'espace $L^2(dr)$. On en d\'eduit que la composante de $\vkn$ selon $y_{-}$ est nulle et donc que $\vkn$ v\'erifie:
\begin{equation}
\label{asymptotique}
    \left(\frac{d}{dr}\right)^j \vkn(r) \underset{r \rightarrow 0}{\lesssim} 
    r^{1-j + 
    \sigma_k},\quad j=0,1,2
\end{equation}
Cette majoration n'est uniforme ni en $n$, ni en $k$, mais elle justifie toutes les int\'egrations par parties qui suivent gr\^ace \`a la remarque \'el\'ementaire suivante:
\begin{remarque}
\label{rem IPP}
Soient $F$ et $G$ deux fonctions r\'eguli\`eres de la variable $r>0$, nulles \`a l'infini, telles que pour des r\'eelles $\sigma$ et $\tau$ v\'erifiant $\sigma+\tau>0$:
\begin{equation*}
 \left(\frac{d}{dr}\right)^j F(r)  \underset{r\rightarrow 0}\lesssim r^{\sigma-j},\quad
 \left(\frac{d}{dr}\right)^j G(r) \underset{r\rightarrow 0} \lesssim r^{\tau-j},\quad j=0,1.
\end{equation*}
Alors:
$$ \int_{0}^{+\infty} \frac{d F}{d r} G\,dr=-\int_{0}^{+\infty} 
F \frac{d G}{d r}\, dr$$
\end{remarque}
Notons $\wkn=e^{i \frac rh}\vkn$ et posons:
\begin{gather*}
\beta_k^2(n)\egaldef\int_{0}^l |\wkn|^2 \d r,\quad \gamma_k^2(n)\egaldef\int_0^l
h^2_n\left|\wkn'\right|^2 \d r,\quad
\kappa_k^2(n)\egaldef\int_0^l |\gkn|^2 \d r\\
M_k^2(n)\egaldef \gamma_k^2(n)+(1+h^2\nu_k^2)\beta_k^2(n)+\kappa_k^2(n).
\end{gather*} 
Les suite $u_n$, $h\nabla u_n$, et $\frac{h}{r} u_n$ sont born\'ees dans $L^2_{\loc}$, et $f_n$ tend vers $0$ dans $L^2$. On a donc:
\begin{equation}
\label{kappak et Mk}
\sum_{k=0}^{+\infty} \kappa_k^2(n) \underset{n\rightarrow +\infty}{\longrightarrow} 0, \quad \sum_{k=0}^{+\infty} M_k^2(n) \underset{n\rightarrow+\infty}{=}O(1).
\end{equation}
L'\'equation v\'erifi\'ee par $\wkn$ pour $r\leq l$ s'\'ecrit:
\begin{equation}
\label{eq wk}
\begin{gathered}
T_k \wkn=he^{i r/h}\gkn\\
T_k\overset{def}{=} e^{ir/h} S_k e^{-ir/h} 
=-h^2 \frac{d^2}{dr^2} +h^2 \frac{b_k}{r^2}+h^2
V+ih\left(\alpha+2\frac{d}{dr}\right)\\
b_k\egaldef \nu_k^2+\frac{d^2-4d+3}{4}
\end{gathered}
\end{equation}
On a:
\begin{equation*}
\nu_k\geq\tilde{\nu}\Rightarrow \sigma_k \geq \tilde{\sigma},\quad
\tilde{\sigma}\egaldef \sqrt{\left(d/2-1 \right)^2 +a+\tilde{\nu}^2}.
\end{equation*}
On se donne une fonction positive $\varphi$ de $C^{\infty}_0([0,l[)$ valant
$1$ pr\`es de $0$. Le r\'eel $t\in]0,1[$ \'etant fix\'e, on choisit $\tilde{\nu}$ tel que:
\begin{equation}
 \label{conditions t}
 \quad 2\tilde{\sigma}-t>0.
\end{equation}
On commence par montrer que si $\tilde{\nu}$ est assez grand:
\begin{equation}
\label{ineg w 1}
\int h_n^2 |\wkn'|^2 r^{-t} \varphi \d r+\int h_n^2|\wkn|^2 r^{-2-t} \varphi \d r=O\left( M^2_k(n)\right),
\end{equation}
o\`u on a not\'e:
$$X_k(n)=O\left( Y_k(n)\right) \iff \exists C>0,\quad\forall k,\;\nu_k\geq \tilde{\nu},\quad \forall n,\quad |X_k(n)|\leq C |Y_k(n)|.$$
La constante $C$ d\'epend donc \'eventuellement de $t$ mais ni de $n$, ni
de $k$, moyennant la condition: $\nu_k\geq\tilde{\nu}$.
On a (dans tous les calculs qui suivent on omet les indices $n$ pour
all\'eger les notations):
\begin{multline}
\label{calcul 1}
\overbrace{\Re\int T_k \wk \wkbar' r^{1-t}\varphi \d r}^{I}=\overbrace{-\Re h^2
  \int \wk'' \wkbar' r^{1-t} \varphi \d r}^{I_1}\\+\underbrace{\Re
  h^2 b_k \int \wk \wkbar' r^{-1-t}\varphi \d r}_{I_2}+
\underbrace{\Re
  h^2 \int V \wk \wkbar' r^{1-t}\varphi \d r}_{I_3}+
\underbrace{\Im h\alpha \int \wk \wkbar'r^{1-t} \varphi \d r}_{I_4}
\end{multline}
Remarquons que par les estimations (\ref{asymptotique}), et l'hypoth\`ese (\ref{conditions t}) sur le param\`etre $t$, toutes ces int\'egrales sont absolument convergentes. On a:
\begin{gather*}
I_1=-\frac{1}{2} h^2 \int \frac{d}{dr} |\wk'|^2 r^{1-t} \varphi \d r\\
\notag
|\wk'|^2\lesssim r^{2\sigma_k-1},\quad \frac{d}{dr} |\wk'|^2\lesssim r^{2\sigma_k-2},\quad r\rightarrow 0.
\end{gather*}
D'apr\`es la remarque \ref{rem IPP} et la condition $2\tilde{\sigma}-t>0$, on peut int\'egrer par parties:
\begin{equation}
\label{calcul 2}
I_1=\underbrace{\frac {(1-t)}{2} h^2\int |\wk'|^2 r^{-t} \varphi \d r}_{I_{1a}}+ \underbrace{\frac 12 h^2 \int |\wk'|^2 r^{1-t} \varphi' \d r}_{O(\gamma_k^2(n))}.
\end{equation}
Pour majorer le dernier terme on a utilis\'e que la d\'eriv\'ee de $\varphi$ est nulle pr\`es de $0$. En raisonnant de la m\^eme mani\`ere pour justifier les int\'egrations par parties, on obtient:
\begin{equation}
\label{calcul 3}
\begin{aligned}
I_2  = & \frac{h^2}{2} b_k \int \frac{d}{dr} |\wk|^2 r^{-1-t} \varphi \d r\\
 I_2= & \underbrace{\frac{(1+t)}{2} h^2 b_k \int |\wk|^2 r^{-2-t} \varphi \d r}_{I_{2a}}-\underbrace{\frac{h^2}{2}b_k \int |\wk|^2 r^{-1-t} \varphi' \d r}_{O(h^2(1+\nu_{k}^2) \beta^2_k)}.
\end{aligned} 
\end{equation}
Par la majoration (\ref{HypV3}) sur $V$, on a:\
\begin{align}
\notag
|I_3|\leq & h^2\int \frac{C_V}{r^2} |\wk||\wk'| r^{1-t}\varphi \d r\\
\notag
\leq & C_V \left\{\int \frac {h^2}{r^2} |\wk|^2 r^{-t} \varphi \d r \right\}^{1/2}
\left\{\int h^2 |\wk'|^2 r^{-t}\varphi \d r \right\}^{1/2}\\
\label{calcul 4}
|I_3|\leq & \frac{C_V}{2}\left\{\frac{1}{\eps} \int \frac {h^2}{r^2} |\wk|^2 r^{-t}
  \varphi \d
  r+\eps \int h^2 | \wk'|^2 r^{-t} \varphi \d r\right\}
\end{align}
Enfin, on a les majorations simples:
\begin{equation}
\label{majorations}
\begin{gathered}
I_4=O(\alpha \gamma_k \beta_k),\\
I=\int T_k \wk \wkbar' r^{1-t}\varphi \d r=\int e^{i\frac{r}{h}} h \gk
\wkbar' r^{1-t}\varphi \d r =O( \kappa_k \gamma_k), 
\end{gathered}
\end{equation}
Prenons $\eps$ assez petit pour que $\frac{1-t}{2}>\frac{C_V\eps}{2}$, ce
qui est possible car $t<1$, puis choisissons $\tilde{\nu}$ tel que:
$$ \nu_k\geq \tilde{\nu}\Rightarrow b_k> \frac{C_V}{2\eps}. $$
Les deux termes principaux de (\ref{calcul 2}) et (\ref{calcul 3}),
$I_{1a}$ et $I_{2a}$  dominent
donc $I_3$ et on obtient, par (\ref{calcul 1}), (\ref{calcul 2}),
(\ref{calcul 3}), (\ref{calcul 4}), (\ref{majorations}), l'in\'egalit\'e
(\ref{ineg w 1}).
 On va maintenant en d\'eduire, en utilisant \`a nouveau l'\'equation
 (\ref{eq wk}) sur $\wk$ la conclusion du lemme \ref{le.lem.ug.en.0}. On a:
\begin{multline}
\label{egalite.J}
J\egaldef \Im \int T_k \wk \wkbar r^{1-t} \varphi \d r\\= \underbrace{-h^2 \Im \int \wk'' \wkbar r^{1-t} \varphi \d r}_{J_1}+\underbrace{\alpha h\int |\wk|^2 r^{1-t}\varphi \d r}_{J_2}+\underbrace{2h \Re\int \wk' \wkbar r^{1-t} \varphi \d r}_{J_3}.
\end{multline}
En remarquant (les int\'egrations par parties se justifient comme pr\'ec\'edemment):
\begin{align*}
J=&O\left(h\kappa_k \beta_k\right)\\
J_1=&2 \Im(1-t)h^2\int \wk' \wkbar r^{-t} \varphi \d r+O(h\beta_k\gamma_k)\\
|J_1|\leq& C h \left\{ \int h^2 |\wk'|^2 r^{-t} \varphi \d r\right\}^{1/2}\left\{ \int |\wk|^2 r^{-t}\varphi \d r\right\}^{1/2} +O(h\beta_k \gamma_k)\\
J_2=&O(h\beta_k^2)\\
J_3=&h\int \frac{d}{dr} |\wk|^2 r^{1-t}\varphi \d r\\ =&(t-1)h\int |\wk|^2r^{-t} \varphi \d r+O(h\beta_{k}^2) 
\end{align*}
En utilisant (\ref{ineg w 1}) pour majorer $J_1$, on d\'eduit de (\ref{egalite.J}):
$$ \int |\wkn|^2 r^{-t} \d r=O(M_k^2(n)).$$
En repassant \`a $\ukn=e^{-ir/h} r^{(d-1)/2} \wkn$, puis en sommant sur tous les $k$
tels que $\nu_k\geq\tilde{\nu}$, et en utilisant la majoration (\ref{kappak et Mk}) de la somme des $M_k^{2}$, on obtient exactement
(\ref{lem.ug.en.0}).
\end{proof} 

\begin{lem}
\label{lem.equivalence}
Soient: 
$l>0$, $q_1,q_2 \in L^{\infty}_{\loc}\left(]0,l],\RR\right)$ tels que:
\begin{gather}
\label{Hyp.q1}
q_1(r)\geq \frac{b}{r^2},\quad b>-\frac{1}{4}\\
\label{Hyp.q}
\quad |q_j(r)|\leq \frac{C}{r^2},\quad j=1,2. 
\end{gather}
Alors l'\'equation diff\'erentielle:
\begin{equation}
\label{eq.diff}
 y''(r)=\big(q_1(r)+iq_2(r)\big) y(r),\quad r\in]0,l],
\end{equation}
admet une base de solutions $(y_+,y_-)$ telle que:
\begin{align}
\label{condition.y_+}
   \left(\frac{d}{dr}\right)^j y_{+}(r) &\underset{r\rightarrow 0}{\lesssim} 
    r^{1/2+\sqrt B-j},\quad j=0,1,2\\
\label{condition.y_-}
    y_{-}(r) &\underset{r\rightarrow 0}{\gtrsim}  
    r^{1/2-\sqrt B}
\end{align}
\end{lem}
\begin{proof}
Posons:
$$ r=e^{-s}, \quad y(r)=e^{-s/2}z(s)=r^{1/2} z(-\log r).$$
L'\'equation diff\'erentielle (\ref{eq.diff}) et les hypoth\`eses (\ref{Hyp.q1}) et (\ref{Hyp.q}) s'\'ecrivent:
\begin{align}
\label{eq.diff'}
\tag{\ref{eq.diff}'}
z''=\underbrace{\left(1/4+e^{-2s} q_1(e^{-s})\right)}_{Q_1(s)} z+&\underbrace{ie^{-2s} q_2(e^{-s})}_{Q_2(s)}z,\quad s\geq L\egaldef -\log l\\
\label{Hyp.q1'}
\tag{\ref{Hyp.q1}'}
Q_2(s) \geq& B\egaldef b+1/4 > 0\\
\label{Hyp.q'}
\tag{\ref{Hyp.q}'}
|Q_j(s)|\leq& C,\quad j=1,2. 
\end{align}
On s'inspire de \cite[Chap. II Ex. 14]{Tosel}. Soit $z$ une solution de (\ref{Hyp.q1'}) et $Z=|z|^2$. Alors, en utilisant (\ref{eq.diff'}) et (\ref{Hyp.q1'}), on obtient:
\begin{align}
\notag
Z'(s)= & 2 \Re \left(z(s)\overline{z}'(s)\right)\\
\notag
Z''(s)= & 2|z'(s)|^2+2 Q_1(s)|z(s)|^2\\
\notag
\geq& 2\sqrt{B}\left( \frac{1}{\sqrt{B}} |z'(s)|^2+\sqrt{B} |z(s)|^2 \right)\\
\label{ineg.Z''}
Z''(s) \geq& 2\sqrt{B} Z'(s).
\end{align}
On commence par construire la solution $z_-$ correspondant \`a $y_-$, en choisissant la solution de (\ref{eq.diff}') telle que $ z_-(L)=z_-'(L)=1.$
Il est facile de voir, avec (\ref{ineg.Z''}) que $Z_-$ et $Z'_-$ restent sup\'erieur \`a $1$ et croissantes. Par le lemme de Gronwall, (\ref{ineg.Z''}) implique:
\begin{align}
\notag
Z_-'(s) \geq & e^{2\sqrt{B}(s-L)} Z_-'(L)\\
\notag
Z_-(s_0)\geq & Z_-'(L) \frac{1}{2\sqrt{B}}e^{2\sqrt{B}(s_0-L)}+Z_-(L),\\
\label{prop.z_-}
z_-(s)\underset{s \rightarrow +\infty}{\gtrsim} & e^{\sqrt B s}
\end{align}
ce qui d\'emontre la condition (\ref{condition.y_-}) sur $y_-(r)=r^{1/2} z_-(-\log r)$. On construit ensuite $z_+$, comme la solution de (\ref{eq.diff'}) d\'efinie par:
$$ z_+(s)\egaldef z_-(s)\int_s^{+\infty} \frac{\d \sigma}{z_-^2(\sigma)}.$$
On a:
\begin{align*}
|z_+(s)| \leq & |z_-(s)|\int_s^{+\infty} \frac{\d \sigma}{|z_-(\sigma)|^2}\\
\leq & \int_s^{+\infty} \frac{\d \sigma}{|z_-(\sigma)|},
\end{align*}
car $|z_+|$ est croissante puisque $Z_+$ l'est. Donc, gr\`ace \`a (\ref{prop.z_-}):
\begin{equation}
\label{prop.z_+}
z_+(s) \underset{s\rightarrow +\infty}{\lesssim} e^{-\sqrt{B}s},
\end{equation}
et par l'\'equation (\ref{eq.diff'}) et les hypoth\`eses sur $Q_1$ et $Q_2$ la m\^eme propri\'et\'e est vraie pour la d\'eriv\'ee seconde de $z_+$. Il est bien connu que cela implique aussi (\ref{prop.z_+}) pour la d\'eriv\'ee premi\`ere de $z_+$. Les estimations (\ref{condition.y_+}) sur $y_+(r)=r^{1/2} z_+(-\log(r))$ en d\'ecoulent imm\'ediatement.
\end{proof}

\subsection{Fin de la d\'emonstration dans le cas multi-polaire}
\label{multipolaire}
Dans cette partie, on ach\`eve la preuve du th\'eor\`eme dans le cas $N\geq
2$ en montrant:
\begin{prop}
\label{dern.prop}
Supposons (\ref{HypV0}),..,(\ref{HypV5}).
Soit $p\in \PPP$ et $\widetilde{\PPP}_p$ l'ensemble des p\^oles $p'$, distincts de
$p$, et tels que $\mu$ ne soit pas nulle au voisinage de $p'$. On suppose que
$\mu$ n'est pas nulle pr\`es de $p$. Alors $p$ appartient \`a
l'enveloppe convexe de $\widetilde{\PPP}_p$.
\end{prop}
Consid\'erons $\widetilde{\PPP}$ l'ensemble des p\^oles pr\`es desquels la mesure
n'est pas nulle. Supposons $\widetilde{\PPP}$ non vide. Le bord de l'enveloppe convexe de
$\widetilde{\PPP}$ est alors un
polyg\^one dont les sommets sont des points de $\widetilde{\PPP}$. En appliquant la
proposition \ref{dern.prop} \`a un tel sommet, on obtient une contradiction qui montre que
$\widetilde{\PPP}$ est vide et donc que $\mu$ est nulle, ce qui prouve le th\'eor\`eme \ref{theo1}.
\begin{proof}
On suppose encore $p=0$ et on reprend les notations des deux parties
pr\'ec\'edentes. Puisque d'apr\`es la proposition
\ref{pas.de.concentration}, $\mu$ ne charge pas $0$, la formule (\ref{exp
  de mu}) peut se r\'e\'ecrire:
 \begin{equation}
    \label{exp de mu'}
    \tag{\ref{exp de mu}'}
\begin{gathered}
    \big<\mu_{\gg},a\big>=\big<\mu,a\big>= \sum_{p_{j}\in\widetilde{\PPP}_0} \int_{0}^{+\infty} 
    \left(\lambda_{i}^{+} a(x=tZ_{i},\xi=Z_{i})+\lambda_{i}^{-} 
    a(x=tZ_{i},\xi=-Z_{i})\right)\,dt\\
a\in C_0^{\infty}\left(\{x\in\RR^d,|x|\leq l\}\times\RR^d_{\xi}\right),
\end{gathered}
\end{equation}
o\`u les $\lambda_i^{\pm}$ sont des constantes positives. On commence par
d\'emontrer, pr\`es de $0$, l'\'equivalent du lemme 
\ref{mu a l'infini} d'invariance globale:
\begin{lem}
    Les rayons rentrant en $0$ et sortant de $0$ portent la m\^eme 
masse:
\begin{equation}
\label{def.Lambda}   
\sum_{p_{j}\in \widetilde{\PPP}_0} \lambda_{j}^{-}=\sum_{p_{j}\in \widetilde{\PPP}_0} 
    \lambda_{j}^{+}\defegal \Lambda. 
\end{equation}
\end{lem}
\begin{proof}
On se donne une fonction $\varphi\in \Cio(\{|x| < l\})$, positive, radiale, 
d\'ecroissante en $r$ et valant $1$ pr\`es de $0$. On a:
\begin{equation*}
    S_n\egaldef h_n^{-1}\big( h_n^2 P u_n-u_n,\varphi 
u_n\big)_{L^2}-h_n^{-1}\big( 
    \varphi u_n, h_n^2 P u_n-u_n\big) \underset{n\rightarrow 
+\infty}{\longrightarrow} 
    0.
\end{equation*}
Mais:
\begin{align*}
    S_n=& h_n Q(u_n,\varphi u_n)-Q(\varphi u_n,u_n)\\
    =& h_n \int \left(\nabla u_n.\nabla \varphi\right) \ubar \d x  -h_n \int u_n \nabla     \varphi.\nabla \ubar_n \d x\\
    =& 2i \Im \int h_n\nabla u_n.\nabla \varphi \ubar_n \d x.
\end{align*}
La suite $(u_n)$ \'etant $h_n$-oscillante, il est facile de 
montrer:
$$ 0=\lim_{n\rightarrow +\infty} S_n= 2i \Im \big< \mu,i\xi.\nabla \varphi 
\big>$$
On en d\'eduit, en utilisant l'expression (\ref{exp de mu'}):
\begin{align*}
    0=&\sum_{p_{j}\in \widetilde{\PPP}_0} \int_{0}^{+\infty} \left( 
\lambda_{i}^{+} 
    \frac{d \varphi}{d r}(r=t)-\lambda_{i}^{-} 
    \frac{d \varphi}{d r} (r=t) \right)\,dt\\
    =&-\sum_{p_{j}\in \widetilde{\PPP}_0} \lambda_{i}^{+} +\sum_{p_{j}\in \widetilde{\PPP}_0} 
    \lambda_{i}^{-}.
\end{align*}
\end{proof}
On \'ecrit maintenant une variante du lemme \ref{le.lem.ug.en.0}, qui
utilise l'hypoth\`ese (\ref{HypV4}) faite sur la d\'eriv\'ee de $V$: 
\begin{lem}
\label{le.lem.important}
Soit $c_d=\frac{(d-1)(d-3)}{4}$. Alors:
\begin{equation}
\label{lem.important}
\left(c_d+\tilde{\nu}-\frac{C_V'}{2}\right) \limsup_{n\rightarrow +\infty}
\int \frac{h_n^2}{|x|^3} |u_{\gg n}|^2 \d x \leq 2 \Lambda.
\end{equation}
\end{lem}
\begin{remarque}
Le lemme \ref{le.lem.important} est similaire \`a la majoration interm\'ediaire
(\ref{ineg w 1}) du lemme \ref{le.lem.ug.en.0}, avec $t=1$, mais on ne peut pas
d\'emontrer la majoration correspondante de la d\'eriv\'ee radiale de
$u_{\gg}$, (\`a cause de la constante $1-t$ dans (\ref{calcul 2}) qui
s'annule lorsque $t=1$):
$$ \limsup_{n\rightarrow +\infty}
\int \frac{h_n^2}{|x|} \left|\frac{\partial u_{\gg n}}{\partial r}\right|^2 \d x <
+\infty.$$
De m\^eme, on ne peut malheureusement pas aboutir par cette m\'ethode \`a
la conclusion (\ref{lem.ug.en.0}) du lemme \ref{le.lem.ug.en.0} dans le cas
$t=1$, c'est \`a dire montrer:
$$ \limsup_{n\rightarrow +\infty} \int |x|^{-1} \left|u_{\gg n}\right|^2\d x
<+\infty.$$
 Notons que cette derni\`ere propri\'et\'e, avec la formule (\ref{exp de mu'})
et la non-int\'egrabilit\'e en $0$ de l'application $r\mapsto 1/r$,
impliquerait directement la nullit\'e de $\mu$ pr\`es de $0$.
\end{remarque}
\begin{remarque}
Le lemme \ref{le.lem.important} correspond \`a un gain d'une demi puissance
de $r$ par rapport \`a la borne naturelle:
$$ \frac{h_n}{r} u_n=O(1) \text{ dans } L^2_{\loc}. $$
Dans la suite, la fonction $h_n^{3/2} r^{-1} u_n$, jouera, pr\`es $0$, un
r\^ole similaire \`a celui que jouerait une trace dans un probl\`eme au
bord. On peut comparer le gain de $r^{-1/2}$ donn\'e par
(\ref{lem.important}) au gain d'une demi-d\'eriv\'ee par rapport au
th\'eor\`emes de traces standard obtenu sur les traces d'un
probl\`eme au bord  avec des ``bonnes'' conditions au bord.
\end{remarque}
\begin{proof}
On se donne $k$ tel que $\nu_k\geq \tilde{\nu}$, et on reprend la d\'emonstration du lemme \ref{le.lem.ug.en.0}, mais avec $t=1$.
La formule (\ref{calcul 1}) est encore valable et s'\'ecrit:
\begin{multline}
\tag{\ref{calcul 1}'}
\label{calcul.1'}
\overbrace{\Re\int T_k \wk \wkbar' \varphi \d r}^{I}=\overbrace{-\Re h^2
  \int \wk'' \wkbar'  \varphi \d r}^{I_1}\\+\underbrace{\Re
  h^2 b_k \int \wk \wkbar' r^{-2}\varphi \d r}_{I_2}+
\underbrace{\Re
  h^2 \int V \wk \wkbar' \varphi \d r}_{I_3}+
\underbrace{\Im h\alpha \int \wk \wkbar' \varphi \d r}_{I_4}.
\end{multline}
Toutes le int\'egrales
\'ecrites sont absolument convergentes gr\^ace aux estimations
(\ref{asymptotique}). Par des int\'egrations par parties \'el\'ementaires,
que l'on justifie encore par la remarque \ref{rem IPP}, on calcule:
\begin{eqnarray}
\label{calcul 1b}
I_1(k,n)&=&\frac{1}{2}\int h^2 |\wk'|^2 \varphi'(r) \d r\\
\label{calul 2b}
I_2(k,n)&=&\underbrace{\int h^2 b_k |\wk|^2 r^{-3} \varphi \d r}_{I_{2a}}+\underbrace{\int -\frac{h^2\nu_k^2}{r^2} |\wk|^2 \varphi'(r) \d r }_{I_{2b}}
+ O(h^2 \beta_k^2).\\
\notag
I_3(k,n)&=& -\frac{h^2}{2} \int \frac{d V}{dr} |\wk|^2 \varphi \d r
+O(h^2\beta_k^2)\\
\label{calul 3b}
|I_3(k,n)|&\leq& \frac{C'_V}{2} \int \frac{h^2}{r^3} |\wk|^2\varphi \d r+
O(h^2\beta_k^2)\\
\label{majoration b}
I_4(k,n)&=& O(\alpha \gamma_k \beta_k).
\end{eqnarray}
Notons: 
$$ \III_j(n)\egaldef \sum_{\nu_k\geq \tilde{\nu}} I_j(k,n),\quad
\III_{2a,b}(n)\egaldef \sum_{\nu_k\geq \tilde{\nu}} I_{2a,b}(k,n).$$
En sommant par rapport \`a l'indice $k$, on a facilement les estimations:
\begin{align}
\label{III2a}
\III_{2a}(n)&\geq (\tilde{\nu}+c_d)\int
\frac{h_n^2}{|x|^3}\left|u_{\gg n}\right|^2 \varphi \d x\\
\label{III3}
|\III_{3}(n)|&\leq \frac{C'_V}{2}\int \frac{h_n^2}{|x|^3} |u_{\gg n}|^2\varphi \d
x+o(1),\quad n\rightarrow +\infty\\
\label{III4}
\III_{4}(n) & \underset{n\rightarrow +\infty}{\longrightarrow} 0.
\end{align}
Il reste \`a \'etudier $\III_1$ et $\III_{2b}$, dont on peut calculer la
limite quand $n$ tend vers $+\infty$ en utilisant la formule (\ref{exp de
  mu'}). On a:
$$ e^{-ir/h} r^{-\frac{d-1}{2}} h \wk'=e^{-ir/h} r^{-\frac{d-1}{2}}
h\frac{d}{dr}\left( e^{ir/h} r^{\frac{d-1}{2}}\uk \right)=h\uk'+i \uk+
h\frac{d-1}{2 r} \uk.$$
En sommant par rapport \`a $k$, on obtient:
\begin{equation}
\label{egalite.I1}
\III_1(n)=\frac{1}{2}\int \left| h\frac{\partial u_{\gg}}{\partial r} +i
  u_{\gg}\right|^2 \varphi'(|x|) \d x+o(1),\quad n\rightarrow +\infty.
\end{equation}
Posons: 
\begin{eqnarray*}
A\egaldef h\frac{\partial}{\partial r}+i,&\quad
 \tilde{A}\egaldef\tilde{\varphi}(x)A\psi(hD),\\
\tilde{\varphi}\in C_0^{\infty}(\RR^d),&\quad \tilde{\varphi}=1 \text{ sur
}\supp \varphi',\\
\psi \in C_0^{\infty}(\RR^d),&\quad \psi(\xi)=1 \text{ sur } \{|\xi|=1\}
\end{eqnarray*}
Les suites $(u_{\gg n})$ et $(h_n \nabla u_{\gg n})$ \'etant
$h_n$-oscillantes sur le support de $\varphi'$, et le support de la mesure
$\mu_{\gg}$ \'etant inclus dans l'ensemble $\{ |\xi|=1\}$, on montre
facilement les convergences $L^2$:
$$ \lim_{n\rightarrow +\infty} \indic_{\supp \varphi'}\left(\psi(h_n D)u_{\gg n}\right)=0,\quad
 \lim_{n\rightarrow +\infty} \indic_{\supp \varphi'}\left(\psi(h_n
  D)\nabla u_{\gg n}\right)=0. $$
Et donc:
$$ \lim_{n\rightarrow 0} \left( \tilde{A}-A \right) u_{\gg n}=0 \text{
  dans }L^2(\supp \varphi').$$
On en d\'eduit:
\begin{align}
\notag
 \III_1(n) =& \frac{1}{2} \left( \tilde{A}^* \varphi'(|x|) \tilde{A}
  u_{\gg n},\tilde u_{\gg n}\right)_{L^2}
+\underset{n\rightarrow+\infty}{o(1)}\\
\notag
\lim_{n\rightarrow \infty} \III_1(n)=&\left\langle \mu, \frac{1}{2} \psi(\xi)
\left| i\xi.\frac{x}{|x|} +i\right|^2 \varphi'(|x|)\right\rangle\\
\notag
=& \sum_{p_j \in \widetilde{\PPP}_0} \frac{1}{2} \int_0^{+\infty} \lambda_i^+ 2^2
\varphi'(t) \d t\\
\label{III1}
\lim_{n\rightarrow +\infty} \III_1(n)=& 2\Lambda,  
\end{align}
o\`u aux deux derni\`eres lignes, on a utilis\'e la formule
(\ref{exp de mu'}), la d\'efinition (\ref{def.Lambda}) de $\Lambda$ et le
fait que $\varphi(0)=1$.
D'autre part:
\begin{align*}
\III_{2b}(n)=\sum_{\nu_k\geq \tilde{\nu}} -\frac{h_n^2}{2} \nu_k^2 \int |\wkn|^2 r^{-2}
\varphi'(r) \d r&=\sum_{\nu_k\geq \tilde{\nu}} -\frac{h_n^2}{2}\nu^2_k\int
\frac{1}{|x|^2}|\ukn|^2\varphi'(|x|) \d x\\
&= \int -\frac{h_n^2}{2} |\nabla_T u_{\gg n}|^2 \varphi'(|x|) \d x 
\end{align*}
o\`u $\nabla_T$ d\'esigne le gradient tangentiel:
$$\nabla_T U=\nabla U-\left(\frac{\partial U}{\partial r}\right) \frac{x}{|x|}.$$ 
La limite de ce terme, quand $n$ tend vers
l'infini, est donc, par un raisonnement analogue au pr\'ec\'edent:
\begin{equation}
\label{III2b}
\lim_{n\rightarrow +\infty} \III_{2b}(n)=
\left<\mu,\frac{1}{2}\varphi'(|x|)
  \left|\xi-(x.\xi)\frac{x}{|x|^2}\right|^2\right>,
\end{equation}
qui vaut $0$ car $\mu_{\gg}$ est nulle lorsque $\xi$ n'est pas paral\`elle \`a $x$.
En remarquant:
$$ \III(n)\egaldef\sum_{j=1}^{4} \III_j(n)=o(1),\quad n\rightarrow
+\infty,$$
 et en  additionnant (\ref{III2a}), (\ref{III3}), (\ref{III4}),
 (\ref{III1}), (\ref{III2b}) on obtient (\ref{lem.important}).
\end{proof}
Introduisons maintenant le r\'eel positif $l_{\gg}$, 
d\'efini par:
$$ l_{\gg}\egaldef\limsup_{n\rightarrow +\infty}\int_{r\leq 1} 
\frac{h_n^2}{|x|^3} |u_{\gg n}|^2\d x.$$
D'apr\`es le lemme pr\'ec\'edent, on a, d\`es que $\tilde{\nu}$ st assez grand:
\begin{equation}
    \label{inegW}
    0\leq l_{\gg}\leq \frac{2 \Lambda}{c_{d}+\tilde{\nu}^2-\frac{C'_V}{2}}
\end{equation}
Les $\lambda_{j}^{+}$ et les $\lambda_{j}^{-}$ vont d\'esormais jouer 
des r\^oles sym\'etriques. La mesure $\mu$ \'etant non nulle au voisinage
de $0$, on sait que $\Lambda$ est strictement positif et on peut donc poser:
\begin{gather*}
 2 \Lambda t_{j}\egaldef 
\lambda_{j}^{+}+\lambda_{j}^{-},\quad 
\sum_{p_j \in \widetilde{\PPP}_0} t_{j}=1,\quad 0\leq t_{j}\leq 1\\
Z\egaldef\sum_{p_j \in \widetilde{\PPP}_0} t_j Z_j,
\end{gather*}
Si l'on suppose que $0$ n'appartient pas \`a l'enveloppe convexe de
$\widetilde{\PPP}_0$, le vecteur $Z$ (qui ne d\'epend pas de $\tilde{\nu}$
d'apr\`es la proposition (\ref{mup.nulle})), est non nul. Notons que si
$\mu$ \'etait invariante par le flot hamiltonien, y compris le long des
rayons passant par le p\^ole, on aurait $Z=0$. Cette invariance se
d\'emontre habituellement en calculant le commutateur $C=[-h^2 P,A]$, o\`u $A$ est un op\'erateur semi-classique $a(x,hD)$, $a\in C^{\infty}_0$. Dans notre cas, $C$ ne peut pas s'estimer, car le potentiel $V$ commute tr\`es
mal avec de tels op\'erateurs. Pour contourner ce probl\`eme, on remplace
$A$ par un op\'erateur diff\'erentiel de degr\'e $1$ dont on peut exprimer
exactement le commutateur avec $V$ en fonction de la d\'eriv\'ee de $V$. L'hypoth\`ese
(\ref{HypV4}), et l'estimation (\ref{inegW}) sur $l_{\gg}$ permettent alors
de ma\^itriser suffisamment $C$ pour montrer que $Z=0$. Ce r\'esultat est
plus faible que le  th\'eor\`eme d'invariance pr\`es du p\^ole, mais suffit
pour d\'emontrer la proposition \ref{dern.prop}.
\begin{lem}
\label{lem.ineg.Z.1}
\begin{equation}
\label{ineg.Z.l}
4\Lambda|Z| \leq \sqrt{d} C'_V l_{\gg}
\end{equation}
\end{lem}
\begin{proof}
Remarquons d'abord que $P$, qui est un op\'erateur born\'e de son domaine $D(P)$ dans $L^2$, d\'efinit par dualit\'e un op\'erateur $\tilde{P}$ born\'e de l'espace $L^2$ dans $D(P)^*$:
$$ (\widetilde{P} F,G)_{D(P)^*,D(P)}=(F,PG)_{L^2},\quad F\in L^2,\quad G\in D(P).$$
De plus, $D(P)$ contient l'espace de fonctions test $C_0^{\infty}(\RR^d\backslash \PPP)$, donc $D(P)^*$ est un sous-espace de l'espace de distributions $\DDD'(\RR^d\backslash \PPP)$, et, d'apr\`es la proposition \ref{D(P)}:
$$ \widetilde{P} F=-\Delta F+ V F \text{ dans } \DDD'\left(\RR^d\backslash \PPP\right).$$
Soit $\Phi=(\Phi_{1},..,\Phi_{d})$ un champ de vecteurs $C^{\infty}$ \`a
support compact dans un petit voisinage de $x=0$ sur $\RR^d$, qui pr\`es de $0$ est constant et \'egal au vecteur $(1,..,1)$. Soit $j$ un entier compris entre $1$ et $d$. On sait que la suite $\left(h_n\Phi_j \partial_{x_j} u_{\gg n}\right)_n$ est born\'ee dans $L^2$. On en d\'eduit:
\begin{align*}
\left( \left(h_n^2 \widetilde{P}-1\right) \Phi_j \partial_{x_j} u_{\gg
    n},u_{\gg n}\right)=& \left(\Phi_j \partial_{x_j} u_{\gg n},\left(h_n^2
    P-1\right) u_{\gg n}\right)\\
=&  o(1), \quad n\rightarrow +\infty. 
\end{align*}
De m\^eme,
\begin{align*}
\left( \Phi_j \partial_{x_j} \left(h_n^2 P-1\right) u_{\gg
    n},u_{\gg n}\right)=& -\left(\left(h_n^2
    P-1\right) u_{\gg n},\partial_{x_j} \left(\Phi_j u_{\gg n}\right)\right)\\
=&  o(1), \quad n\rightarrow +\infty. 
\end{align*}
Et donc:
\begin{equation}
\label{calcul.I.1}
\CCC \egaldef \left(\Phi_j \partial_{x_j} h^2 Pu_{\gg}-h^2 \widetilde{P} \Phi_j \partial_{x_j} u_{\gg},u_{\gg} \right) \underset{n\rightarrow +\infty}{\rightarrow} 0
\end{equation}
Dans (\ref{calcul.I.1}), le produit scalaire est bien d\'efini au sens de
la dualit\'e entre $D(P)^*$ et $D(P)$. On peut exprimer $\CCC$ en fonction de $\mu$ et d'un terme domin\'e par $l_{\gg}$ \`a l'infini:
\begin{align}
\label{calcul.I.2}
\CCC=&\Big(\left[\Phi_j\partial_{x_j},h^2(-\Delta+V) u_{\gg}\right],u_{\gg}\Big)\\
\label{calcul.I.3}
\Big(\left[\Phi_j\partial_{x_j},-h^2 \Delta\right] u_{\gg},u_{\gg}\Big)=&\underbrace{2h^2\big(\nabla \Phi_{j}.\nabla \partial_{x_{j}}u_{\gg},u_{\gg}\big)}_{\longrightarrow \big<\mu,2\xi.\nabla\Phi_j\xi_j\big>}+o(1) \quad n\rightarrow +\infty\\
\notag
\left|\Big(h^2\left[ \Phi_j
      \partial_{x_j},V\right]u_{\gg},u_{\gg}\Big)\right|=& \left|h^2 \left( \Phi_j (\partial_{x_j} V) u_{\gg},u_{\gg}\right)\right|\\
\label{calcul.I.4}
\leq& C'_V l_{\gg} +o(1) \quad n\rightarrow +\infty.
\end{align}
o\`u \`a la derni\`ere ligne on a utilis\'e l'hypoth\`ese
(\ref{HypV4}) et la d\'efinition de
$l_{\gg}$. Finalement, on remarque que d'apr\`es (\ref{exp de mu'}), et
puisque $\Phi_j(0)=1$,
$$ \big<\mu,2\xi.\nabla\Phi_j\xi_j\big>= -4\Lambda Z.e_j,$$
o\`u $e_j$ est le vecteur unitaire dont la $j$-i\`eme coordonn\'ee vaut $1$. En utilisant (\ref{calcul.I.1}),(\ref{calcul.I.2}),(\ref{calcul.I.3}),(\ref{calcul.I.4}), on obtient:
$$ 4 \Lambda|Z.e_j|\leq C_V'l_{\gg},$$
dont on d\'eduit (\ref{ineg.Z.l}) en sommant les carr\'es.
\end{proof}
En utilisant (\ref{inegW}) et le lemme \ref{lem.ineg.Z.1}, on obtient:
$$ 4\Lambda |Z| \leq \sqrt{d} C'_V
\frac{2\Lambda}{c_d+\tilde{\nu}^2-\frac{C'_V}{2}}.$$
Or $Z$ ne d\'epend pas de $\tilde{\nu}$, que l'on peut prendre aussi grand
que l'on veut. On a donc $Z=0$, ce qui montre que le p\^ole $0$ est un barycentre des directions $Z_j$ correspondant au p\^oles $p_j$ qui appartiennent \`a
$\widetilde{\PPP}_0$. On peut en d\'eduire facilement que $0$ est dans l'enveloppe convexe de $\widetilde{\PPP}_0$, ce qui d\'emontre la proposition \ref{dern.prop}. 
\end{proof}

\section{Une variante du th\'eor\`eme \ref{theo1}}
Dans cette section, on \'ecrit une g\'en\'eralisation partielle du th\'eor\`eme \ref{theo1} qui permet de consid\'erer des potentiels de la forme:
$$\frac{1}{|x-p_j|^2} a_j\left(\frac{x-x_j}{|x-x_j|^2}\right) ,\quad a_j+\left(\frac d2 -1 \right)^2>0$$
pr\`es de chaque p\^ole $p_j$. Nous pr\'ecisons ensuite les modifications \`a apporter \`a la d\'emonstration du th\'eor\`eme \ref{theo1} pour d\'emontrer ce r\'esultat.
\label{sec.theo2}
\begin{theo}
\label{theo2}
Soit $\PPP=\{p_1,..,p_N\}$ un ensemble de p\^oles sur $\RR^d$ et $P=-\Delta+V$, o\`u $V$ est un potentiel {\bf r\'eel} sur $\RR^d$ tel que:
$$V=V_1+V_2,$$
avec:
\begin{itemize}
\item $V_1$ v\'erifiant les hypoth\`eses (\ref{HypV0}),\dots,(\ref{HypV5}) du th\'eor\`eme \ref{theo1},
\item $\displaystyle V_2\in L^{\infty}_{\loc}\left(\RR^d\backslash \PPP\right)$, \`a support compact sur $\RR^d$ et tel que, si $|x-p_j|\leq l$:
\begin{gather}
\label{HypV_2} V_2(x)= \frac{1}{|x-p_j|^2} b_j\left(\frac{x-x_j}{|x-x_j|^2}\right),\quad
b_j\in L^{\infty}(S^{d-1})\\
\label{HypV2'}
\tag{\ref{HypV2}'}
b_j\geq 0,\\
\tag{\ref{HypV4}'}
\label{HypV4'}
\nabla_{\theta} b_j  \in L^{\infty}(S^{d-1}).
\end{gather}
\end{itemize}
Alors la conclusion du th\'eor\`eme \ref{theo1} reste valide: si $\chi \in C^{\infty}_0\left(\RR^d\right),$
\begin{equation}
\label{ineg.theo2}
        \exists \lambda_{0}>0,\; \exists C>0,\; \forall \lambda >\lambda_{0},\; \forall 
        \eps>0,\;\quad ||\chi R_{\lambda\pm i\eps}\chi||_{L^2\rightarrow 
        L^2}\leq \frac{C}{\sqrt{\lambda}}.
    \end{equation}
De plus, comme dans le th\'eor\`eme \ref{theo1}, si $N=1$, (\ref{ineg.theo2}) reste vraie sans les hypoth\`eses (\ref{HypV4}) et (\ref{HypV4'}) faites sur le gradient de $V$.
\end{theo}
\begin{remarque}
Un tel $V$ v\'erifie les hypoth\`eses (\ref{HypV0}),\dots,(\ref{HypV3}) du th\'eor\`eme \ref{theo1}, ce qui permet de d\'efinir comme dans la section \ref{sec2}, l'op\'erateur $P$ qu sens des formes quadratiques.
\end{remarque}
\begin{remarque}
Dans le cas unipolaire, les potentiels de la forme:
$$ V(r,\theta)=\chi(x)\frac{a(\theta)}{r^2},\quad a\in L^{\infty}\left(S^{d-1}\right),\quad  a(\theta)\geq a_0,\quad a_0+(d/2-1)^2>0, $$
o\`u $\chi$ est une fonction \`a support compact valant $1$ pr\`es de $0$, rentrent dans le cadre du th\'eor\`eme \ref{theo2}. Contrairement aux r\'esultats connus jusqu'alors, on ne fait aucune hypoth\`ese dans ce cas sur la d\'eriv\'ee du potentiel pr\`es du p\^ole, ce qui est naturel: on n'a pas besoin d'un telle condition pour que $P$ soit positif. Malheureusement, comme dans le cas radial, les conditions sur $|\nabla V|$ r\'eapparaissent lorsque $N\geq 2$ pour \'eliminer les rayons \'eventuels entre les p\^oles du support de $\mu$.
\end{remarque}
\begin{proof}
On reprend point par point la d\'emonstration du th\'eor\`eme \ref{theo1} de la section \ref{sec3}.\par
La partie  \ref{intro.mu}, qui ne n\'ecessite que les hypoth\`eses (\ref{HypV0}),\dots,(\ref{HypV3}), reste valable.\par
On se place ensuite pr\`es d'un p\^ole $p_{j_0}$ que l'on prend comme origine du rep\`ere. L'\'equation sur $u_n$ s'\'ecrit, en coordonn\'ees polaires au voisinage de ce p\^ole:
$$ -h_n ^2 \frac{\partial^2}{\partial r ^2} u_n -h_n^2 \frac{d-1}{r}\frac{\partial}{\partial r}u_n+h_n^2 V_1(r) u_n -\frac{h_n^2}{r^2} \Delta_{\theta} u_n +\frac{b_j(\theta)}{r^2} u_n=h_n\chi f_n,$$
o\`u on a not\'e $\Delta_{\theta}$ le Laplacien sur la sph\`ere $S^{d-1}$. On consid\`ere alors l'op\'erateur sur $L^2(S^{d-1})$:
$$ H=-\Delta_{\theta} +b_j(\theta), \quad D(H)=H^2\left(S^{d-1}\right).$$
La fonction $b_j$ \'etant positive est born\'ee, il est facile de voir que $H$ est un op\'erateur auto-adjoint positif \`a r\'esolvante compacte. Il admet donc une famille de vecteur propre formant une base hilbertienne de l'espace $L^2\left(S^{d-1}\right)$:
\begin{gather*}
   \estar_{k}\in C^{1}(S^{d-1}),\quad H\estar_k=\nu_{k}^{\star 2} e_{k}\\
    \nuu_{k+1}\geq \nuu_{k}\geq \nuu_{0}\geq 0,\quad 
    \nuu_{k}\underset{k\rightarrow +\infty}{\longrightarrow}+\infty  
\end{gather*}
On d\'ecompose, comme dans la partie \ref{petites} chacune des fonctions $u_n$ et $f_n$ selon ces fonctions propres: 
\begin{equation*}
    u_{n}(x)=\sum_{k\in \NN} \uknu(r)\estar_k(\theta),\quad 
    f_{n}(x)=\sum_{k\in\NN}  \fknu(r)\estar_k(\theta).
\end{equation*}
Puis, ayant fix\'e un entier naturel $\tilde{\nu}$, on regroupe \`a nouveaus selon les grandes et les petites valeurs propres:
\begin{equation*}
\begin{gathered}
    u_{n}=u^{\star}_{\pp n}+u^{\star}_{\gg n}, \quad u^{\star}_{\pp n}\overset{def}{=} 
\sum_{\nuu_{k}\leq 
    \tilde{\nu}} \uknu\estar_k,\quad 
     u^{\star}_{\gg n}\overset{def}{=} \sum_{\nuu_{k}>
    \tilde{\nu}} \uknu\estar_k\\
    f_n=f^{\star}_{\pp n}+f^{\star}_{\gg n}, \quad f^{\star}_{\pp n}\overset{def}{=} 
\sum_{\nuu_{k}\leq 
    \tilde{\nu}} \fknu\estar_k,\quad 
     f^{\star}_{\gg n}\overset{def}{=} \sum_{\nu_{k}>
    \tilde{\nu}} \fknu\estar_k.
\end{gathered}
\end{equation*}
Les parties \ref{petites}, \ref{concentration}, \ref{multipolaire} de la d\'emonstration fonctionne alors parfaitement en rempla\c cant les $e_k$ par les $\estar_k$, et les d\'ecompositions en harmoniques sph\'eriques par les d\'ecompositions correspondantes associ\'ees aux $\estar_k$. Remarquons que l'\'equation sur $\uk^{\star}$ pr\`es du p\^ole s'\'ecrit:
$$ -h^2 \uk^{\star}{}'' -h_n^2 \frac{d-1}{r}\uk^{\star}{}'+h_n^2 \left(V_1(r)+\frac{\nu_k^{\star}}{r^2}\right) \uk^{\star}-(1-i\alpha_n)\uk^{\star}=h_n \fk^{\star}.$$   
Elle est donc identique \`a l'\'equation (\ref{eq.uk}) v\'erifi\'ee par $\uk$ dans la d\'emonstration du th\'eor\`eme \ref{theo1}, et en sommant ces \'equations, on obtient pr\`es du p\^oles les bonnes \'equations sur $u_{\gg}^{\star}$ et $u_{\pp}^{\star}$:
\begin{equation}
\label{eq.pg.star}
 h^2P u_{\gg}^{\star}- (1-ih \alpha)u_{\gg}=h f_{\gg}^{\star},\quad h^2P u_{\pp}^{\star}- (1-ih \alpha)u_{\pp}=h f_{\pp}^{\star}.
\end{equation}
Notons que dans les parties \ref{petites} et \ref{concentration}, les calculs se font sur chacune des harmoniques sph\'eriques et que dans \ref{multipolaire}, on utilise seulement l'\'equation sur $u_{\gg}$ identique \`a (\ref{eq.pg.star}), et les hypoth\`eses (\ref{HypV3}) et (\ref{HypV4}) sur $V$, qui sont encore valables ici gr\^ace aux hypoth\`eses (\ref{HypV3}) et (\ref{HypV4}) sur $V_1$, et (\ref{HypV_2}), (\ref{HypV4'}) sur $V_2$. On n'utilise jamais directement, dans \ref{multipolaire}, la forme radiale de $V$ pr\`es du p\^ole. 
\end{proof}


\begin{thebibliography}{9}

\bibitem{ABG} Werner O. Amrein, Anne Boutet de Monvel, Vladimir Georgescu.
\newblock $C^0$-groups, commutator methods and spectral theory of N-body
Hamiltonians. (1996)
\newblock {\em Progress in mathematics. Basel Boston MA Berlin Birkhäuser}

\bibitem{BAK} Ben-Artzi, Matania; Klainerman, Sergiu.
\newblock Decay and regularity for the Schrödinger equation.
\newblock {\em J. Anal. Math. } Vol. 58, 25-37 (1992)

\bibitem{NB1} Nicolas Burq. 
\newblock Semi-classical estimates for 
    the resolvent in non trapping geometries. {\em Int. Math. Res. Not. 2002}, no. 5, 221--241.

\bibitem{BGT} Nicolas Burq, Patrick G\'erard, Nikolai Tzvetkov. 
\newblock On Non-Linear Schr\"odinger Equations in Exterior Domains. Preprint, universit\'e Paris Sud.

    \bibitem{BuLe} Nicolas Burq, Gilles Lebeau. 
\newblock Mesure de d\'efaut de 
    compacit\'e, application au syst\`eme de Lam\'e. {\em Ann. Sci. Ecole Norm. Sup.} (4) Vol. 34, No 6, 817-870 (2001)

    \bibitem{BPSTZ} Nicolas Burq, Fabrice Planchon, John G. Stalker, 
A. 
    Shadi Tahvildar-Zadeh. 
\newblock Strichartz Estimates for the Wave and 
    Schr\"odinger Equations with the Inverse-Square Potential. (2002)

\bibitem{DAP} Piero D'Ancona, Vittoria Pierfelice. 
\newblock On the wave equation with a large rough
potential. Pr\'epublication. (2003)

\bibitem{CEx} Thomas Duyckaerts.
\newblock A Singular Critical Potential For The Schr\"odinger Operator. Prepublication (2003) 
        \bibitem{CFK} Clotilde Fermanian. 
\newblock Equation de la chaleur et mesures 
semi-classiques. Th\`ese de l'universit\'e de Paris-Sud (1995) 

\bibitem{CG} Christian G\'erard, Andr\'e Martinez.
\newblock Principe d'absorption limite pour des op\'erateurs des
Schr\"odinger \`a longue port\'ee. {\em C.R.A.S. S\'erie I,
  math\'ematiques.} 306, No 3, 121-123 (1988)


\bibitem{PG} Patrick G\'erard. 
\newblock Mesures semi-classiques et ondes de Bloch, S\'eminaire EDP de l'Ecole Polytechnique 16 (1990-1991)

\bibitem{Tosel} St\'ephane Gonnord, Nicolas Tosel. Calcul diff\'erentiel. {\em ellipses}


\bibitem{LaPh} P. Lax, R. Phillips. 
\newblock Scattering Theory. 
{\em Pure and Applied Mathematics}, No 26. Academic Press, 2nde edition.

\bibitem{LP} Pierre-Louis Lions et Thierry Paul. 
\newblock Sur les mesures 
de Wigner. {\em Revista Matem\'atica Iberoamericana 9}, n 3, 553-618 (1993)

\bibitem{Mor1} Cathleen S. Morawetz. 
\newblock The Decay for Solutions of the Exterior Initial-Boundary Value Problem for the Wave Equation. {\em Communications on Pure and Applied Mathematics}, Vol. 14, 561-568 (1961)

\bibitem{Mor2} Cathleen S. Morawetz. 
\newblock Decay for Solutions of the Exterior Probem for the Wave Equation. {\em Communications on Pure and Applied Mathematics}, Vol. 27, 229-264 (1975).
\bibitem{Mou} E. Mourre. 
\newblock Absence of Singular Continuous Spectrum for Certain Self-Adjoint Operators. {\em Commun. Math. Phys.} Vol.78, 391-408 (1981)

\bibitem{VP} Vittoria Pierfelice. Decay estimate for the wave equation
  with a small potential. Pr\'epublication. (2003)

\bibitem{PSTZ} Fabrice Planchon, John G. Stalker, A. Shadi Tahvildar-Zadeh. 
\newblock $L^p$ estimates for the wave equation with the inverse-square potential. {\em Discrete Contin. Dyn. Syst.} Vol. 9, No 2, 427-442 (2003)


    \bibitem{RS1} Michael Reed, Barry Simon. 
\newblock Methods of modern 
    mathematical physics. vol I Functionnal analysis. {\em Academic Press}.

    \bibitem{RS2} Michael Reed, Barry Simon. 
\newblock Methods of modern 
    mathematical physics. vol II Fourier analysis, self-adjointness. 
{\em Academic Press}, 1975.

\bibitem{RV1} Alberto Ruiz, Luis Vega. 
\newblock On Local Regularity of Schr\"odinger Equations. {\em Internat. Math. Res. Notices} 1993, 13-27.

\bibitem{RV2} Alberto Ruiz, Luis Vega. 
\newblock Local Regularity of Solutions to Wave Equations with Time-Dependent Potentials. {\em Duke Mathematical Journal.} Vol 76, No 3, 913-940 (1994)
 
\bibitem{VZ} Andr\'as Vasy, Maciej Zworski. 
\newblock Semi-Classical Estimates in Asymptotically Euclidian Scattering. 
{\em Comm. Math. Phys.} Vol. 212, No 1, 205-217 (2000)

\end{thebibliography}
\end{document}